\documentclass[11pt]{amsart}
\usepackage{amsmath}
\usepackage{amsfonts}
\usepackage{amssymb}
\usepackage{amscd}
\usepackage{color}
\input{epsf.sty}
\catcode `\@=11
\def\numberbysection{\@addtoreset{equation}{section}
         \renewcommand{\theequation}{\thesection.\arabic{equation}}}
\numberbysection
\def\subsubsection{\@startsection{subsubsection}{3}%
  \normalparindent{.5\linespacing\@plus.7\linespacing}{-.5em}%
  {\normalfont\bfseries}}
\newlength{\ghost}

\newtheorem{introthm}{Theorem}
\newtheorem{thm}{Theorem}[section]
\newtheorem{lem}[thm]{Lemma}
\newtheorem{prop}[thm]{Proposition}
\newtheorem{cor}[thm]{Corollary}
\newtheorem{claim}[thm]{Claim}

\theoremstyle{definition}

\newtheorem{df}[thm]{Definition}

\newtheorem{rmk}[thm]{Remark}
\newtheorem{nota}[thm]{Notation}
\newtheorem{ex}[thm]{Example}
\newtheorem*{cav}{Caveat}

\newcommand{\set}[2]{\{#1,\dots, #2\}}

\def\In{In}

\def\M{\mathcal{M}}

\def\Sn{\mathbb{S}_n}
\def\Sm{\mathbb{S}_m}
\def\Snn{\mathbb{S}_{n+1}}

\def\Zz{\mathbb{Z}/2\mathbb{Z}}

\def\Z{\mathbb{Z}}

\def\ra{\rightarrow}
\def\del{\partial}
\def\s{\sigma}
\def\D{\Delta}

\def\G{\Gamma}
\def\a{\alpha}
\def\b{\beta}
\def\g{\gamma}

\def\Cact{\mathcal{C}act}
\def\Cacti{\mathcal{C}acti}

\def\pair{\la \;,\;\ra}

\def\del{\partial}

\def\a{\alpha}
\def\ba{\bar{\a}}
\def\b{\beta}
\def\g{\gamma}
\def\bg{\bar{\g}}

\def\D{\Delta}

\def\s{\sigma}

\def\la{\langle}
\def\ra{\rangle}

\def\Sn{\mathbb{S}_n}
\def\Sm{\mathbb{S}_m}
\def\Snn{\mathbb{S}_{n+1}}

\def\comp{Comp}

\def\PP{\mathbb{P}}

\def\lab{\mathrm{lab}}

\def\graphs{\mathcal{G}}
\def\arcgraphs{\overline{\mathcal{G}}}

\def\carcgraphs{\arcgraphs^{e}}

\def\Aang{A^{\angle}}
\def\DA{\mathcal{D}\A}
\def\Darc{\mathcal{DA}rc}
\def\Arc{\mathcal{A}rc}
\def\Arcn{{\mathcal{A}rc_{\#}}}
\def\Darcn{{\mathcal{DA}rc_{\#}}}
\def\Arcno{{\mathcal{A}rc_{\#}^0}}
\def\DArcno{{\mathcal{DA}rc_{\#}^0}}
\def\Tree{\mathcal{T}ree}
\def\Lintree{\mathcal{L}Tree}
\def\Cact{\mathcal{C}act}
\def\Cacti{\mathcal{C}acti}
\def\SCC{\mathcal{S}CC}
\def\Corol{\mathcal{C}orol}
\def\OCArc{\mathcal{C}_o^*(\Arc)}
\def\OCArcn{\mathcal{C}_o^*(\Arc_{\#})}
\def\OCAnarc{\mathcal{C}_o^*(\Anarc)}
\def\OCArcno{\mathcal{C}_o^*(\Arc^0_{\#})}
\def\OCDiarc{\mathcal{C}_o^*(\Diarc)}
\def\GrOCArc{Gr\mathcal{C}_o^*(\Arc))}
\def\GrOCArcno{Gr\mathcal{C}_o^*(\Arc^0_{\#})}

\def\GrOCDiarc{Gr\mathcal{C}_o^*(\Diarc)}

\def\Gr{Gr}
\def\OC{\mathcal{C}_o^*}

\def\Ana{\A^{\angle}}
\def\A{\mathcal{A}}
\def\Agrs{\A_{g,r}^{s}}

\def\Sularc{\mathcal{A}rc^{Sul}}
\def\StSularc{{\mathcal{A}rc^{StSul}}}
\def\Diarc{{\mathcal{A}rc^{i/o}}}
\def\Diarcno{{\mathcal{A}rc^{i/o,0}_{\#}}}
\def\DiA{\A^{i/o}}

\def\Diioarc{\overline{\Arc}^{i\leftrightarrow o}}
\def\LDiioarc{\overline{\L\Arc}^{i\leftrightarrow o}}
\def\Diiooarc{\overline{\Arc}^{i\nleftrightarrow i}}
\def\LDiiooarc{\overline{\L\Arc}^{i\nleftrightarrow i}}
\def\Diiooarci{\overline{\Arc}^{i\nleftrightarrow i}_1}
\def\Diioarci{\overline{\Arc}^{i\leftrightarrow o}_1}
\def\LDiioarci{\overline{\L\Arc}^{i\leftrightarrow o}_1}

\def\iooarc{\mathcal{A}rc^{i \nleftrightarrow i}}
\def\ioarc{\mathcal{A}rc^{i\leftrightarrow o}}
\def\Anarc{\mathcal{A}rc^{\angle}}
\def\Anarcn{\mathcal{A}rc_{\#}^{\angle}}

\def\Digraph{\mathcal{R}ib^{i/o}}
\def\ppdigraph{\mathcal{R}ib^{i \leftrightarrow o}}
\def\StSulchord{\mathcal{C}hord^{StSul}}
\def\Sulchord{\mathcal{C}hord^{Sul}}

\def\Rib{\mathcal{R}ib}
\def\MRib{\mathcal{MR}ib}
\def\PRib{\mathbb{P}\mathcal{R}ib}
\def\Anrib{\mathcal{R}ib^{\angle}}

\def\Loop{\mathcal{L}oop}

\def\mk{\mathrm{mk}}

\def\amark{mk^{\angle}}

\def\val{val}
\def\Ainf{A_{\infty}}

\def\mk{mk}
\def\io{i/o}

\def\In{In}
\def\Out{Out}
\def\punc{pct}

\def\lab{{Lab}}

\def\L{\mathcal{L}}

\def\Loop{\mathcal{L}oop}

\def\Arc{\mathcal{A}rc}

\def\Cyc{\curvearrowright}

\def\Rp{{\mathbb{R}_{>0}}}
\def\PMC{\mathcal{PMC}}

\def\C{Cyc}

\def\Mngn{M_{g,n+1}^{1^{n+1}}}

\newcommand\crl[1]{\vskip 5\lineskip \noindent {\bf Corollary. }{\em #1}
\vskip 5\lineskip}

\begin{document}

\title[Moduli space actions on Hochschild co-chains]
{Moduli space actions on the Hochschild Co-Chains of a Frobenius
algebra I: Cell Operads}

\author
[Ralph M.\ Kaufmann]{Ralph M.\ Kaufmann}
\email{kaufmann@math.uconn.edu}

\address{University of Connecticut, Storrs CT 06269}
\begin{abstract}
This is the first of two papers in which we prove that a cell
model of the moduli space of curves with marked points and tangent
vectors at the marked points acts on the  Hochschild co--chains of
a Frobenius algebra. We also prove that a there is dg--PROP action
of a version of Sullivan Chord diagrams which acts on the
normalized Hochschild co-chains of a Frobenius algebra. These
actions lift to operadic correlation functions on the co--cycles.
 In particular,
the PROP action gives an action on the homology of a loop space of
a compact simply--connected manifold.

In this first part, we set up the topological operads/PROPs and
their cell models. The main theorems of this part are that there
is a cell model operad for the moduli space of genus $g$ curves
with $n$ punctures and a tangent vector at each of these punctures
and that there exists a CW complex whose chains are isomorphic to
a certain type of Sullivan Chord diagrams and that they form a
PROP. Furthermore there exist weak versions of these structures on
the topological level which all lie inside an all encompassing
cyclic (rational) operad.
\end{abstract}

\maketitle


\section*{Introduction}

Recently the operations of cells on the Hochschild complex of
associative algebras have been intensely studied. There are three
main sources for this interest. The first is Deligne's conjecture
which has been solved in various ways
\cite{Maxim,T,MS,MS2,MScosimp,Vor2,KS,MS2,BF,del} together with
its generalization to $A_{\infty}$--algebras and to the cyclic
case using the framed little discs operad \cite{cyclic} in
conjunction with Frobenius algebras. The second source is the
string topology of Chas and Sullivan
\cite{CS,Vor,CJ,C,CG,Chat,Merk,C1,C2}, and lastly the third source
are considerations of $D$-branes in open/closed string theory as
deformations of the Hochschild complex see e.g.\
\cite{KR,KLi1,KLi2}
\footnote{A more extensive discussion of these
links is given in \cite{hoch2}.}. See also \cite{KS2,cost} for
related discussions of these topics in  different settings.

In this paper and its sequel \cite{hoch2} we prove the
following:

\begin{introthm}
\label{introa}
 The moduli space $\Mngn$ of genus $g$ curves with
$n$ punctures and a tangent vector at each of these punctures
has the structure of a rational cyclic operad. This structure
induces a cyclic $dg$ operad structure on a cell model computing the
cohomology of $\Mngn$.

Furthermore the cell level operad operates on the Hochschild
co--chains of a Frobenius algebra. It also yields correlation
functions on the tensor algebra of the co-cycles of a differential
algebra $(A,d)$
 with a cyclically invariant trace $\int: A\to k$
which satisfies $\int da=0$ and whose induced pairing on
$H=H(A,d)$ turns $H$ into a Frobenius algebra.
\end{introthm}

The first part of the theorem together with the graph description
of this moduli space given below can be taken to say that we
define a combinatorial version of conformal field theory.

Furthermore, there is also a PROPic version of this action related
to string topology.

\begin{introthm}{There is a rational topological quasi--PROP which
 is homotopic to a CW
complex whose cellular chains are isomorphic as a free Abelian group
to a certain type of Sullivan Chord diagrams. These chains form a
dg--PROP and hence induce this structure on the Chord diagrams.
Furthermore if $H$ is a Frobenius algebra there is a PROPic action
on the Hochschild co--chains of $H$ that is a dg--action. This
dg-action of a dg-PROP on the dg--algebra of Hochschild co-chains
naturally descends to an action of the homology of the CW-complex on
the Hochschild cohomology of a Frobenius algebra.

Moreover for $(A,d,\int,H)$ as in Theorem \ref{introa} the action
on $H$ is induced by correlation functions on the tensor algebra
of $A$ that yield operadic correlation functions on the tensor
algebra of the co--cycles of $A$ for any $(A,d)$ as above.

Finally, the $BV$ operator, which is given by the action of the
sub--PROP equivalent to the framed little discs operad, acts as in
\cite{cyclic}. Thus the BV operator for the action on the
Hochschild cohomology of $H$ is identified with Connes operator
$B$ under the identification of the Hochschild cohomology of a
Frobenius algebra with its cyclic cohomology of $H$.}
\end{introthm}

The notions  ``quasi'' and ``rational'' denote certain weakenings
of the axioms which are explained below. The definition of the new
notion of operadic correlation functions is contained in the
second part \cite{hoch2}. This notion should be thought of as the
correct definition of a dg--algebra $(A,d)$ over a cyclic operad.
It is the mathematical incarnation of the fact that OPEs in
physics are actually only defined within correlators and only on
BRST closed fields. An immediate consequence
 using Jones' \cite{jones} cyclic description
of the free loop spaces and its $S^1$--action then is

\crl {When taking field coefficients, the above action gives a
$dg$--action of a $dg$--PROP of
 Sullivan Chord diagrams on the $E^1$--term of a spectral sequence
converging to $H_*(LM)$, that is the homology of the loop space of
a simply connected compact manifold, and hence induces operations
on the homology of the loop space.}

Lastly, there is a version of these statements in the case where
$(V,d)$ is a vector space with a differential a pairing $\pair$ that is symmetric
such that $\forall v,w\in V: \la dv, w\ra +\la v,dw\ra =0$ and $\pair$ is   non--degenerate on
$H=H(V,d)$ which is finite dimensional.

\begin{introthm}
The operads and PROPs above also act on the tensor algebra $TV$ of
a triple $(V,d,\pair)$ as specified above and yield operadic
correlation functions for the co-cycles of $TV$.
\end{introthm}

This action is different from the algebra case of Theorems A and B, making the result
interesting in its own right. Furthermore, this action
``descends'' to an action of the stabilized arc operad, which
forms a spectrum \cite{Ribbon}.

The proof of these facts consists mainly of two steps. First
defining the respective topological objects and then defining
their actions.

The first step which is the content of this paper is actually
quite involved, since there are several natural generalizations of
the operads $\Cact$ and $\Cacti$ (see \cite{cact}) yielding the
Gerstenhaber and BV--structures and several types of chain level
models for them. The most notable ones are the $\Arc$ operad of
\cite{KLP} and its subspace $\Arcno$ consisting of quasi--filling
arc families on surfaces without punctures. We prove below that
the space $\Arcno$ is isomorphic to the moduli space of genus $g$
curves with $n$ punctures and a tangent vector at each of these
punctures modded out by the action of $\Rp$ which scales all the
tangent vectors simultaneously. Other generalizations are isomorphic to
different versions of metric Sullivan--Chord diagrams as we explain below.

The reader who is mainly interested in how the action is defined
and wants to forgo the geometric topological and algebraic
topological constructions of the various relevant operads and PROP
can skip ahead to the second part \cite{hoch2}.

In this part, we will systematically extend and augment the operad
structure of $\Arc$ and use these topological structures to induce
operations on the chain level. The results are novel operadic and
PROPic structures associated to various restrictions and extensions
of the moduli space $\Mngn$ on the topological and chain level. The
elements of $\Arc$ or better of its ambient CW complex $\A$ can be
thought of as metric graphs on a surface with boundaries and one
marked point on each boundary, where the set of these points is
equal to  the set of the vertices of the graph. $\Arc$ is then the
subset where all boundaries are ``hit'' that is none of the vertices
has valence $0$. Although $\Arc$ is an operad, $\Arcno$ is only an
operad on a dense subset. The cause of this is that on a codimension
$1$ subset the gluing defined in the ambient $\Arc$ will lead
outside of $\Arcno$. Thus we are lead to consider rational
topological operads, that is operads whose gluings are defined on
dense open subsets. Since the ``bad'' part is of codimension at
least one, this structure, however, descends to a true operad
structure on the chain level for suitable chains. Considering the
dg--structure of the PROP action, we will see in \cite{hoch2}, that
we also have to extend the gluings of the $\Arc$ operad to its
ambient CW complex $\A$ (we give these constructions in \S
\ref{arcsection}).

In order to go into the PROP framework, we need to distinguish the
boundaries as ``in'' and ``out'' boundaries. The resulting space
which accommodates the extra markings is denoted by $\Diarc$. There
are natural additional restrictions which can be imposed on the
graphs. The first is allowing edges only between ``in'' and ``out'';
this space will be called $\ioarc$. Only barring edges from ``in''
to ``in'', we obtain a subspace which will be called $\iooarc$.
These spaces naturally form di--operads. In the dual graph
terminology they correspond to various generalizations of Sullivan
Chord diagrams. We include a discussion of  several versions of
Sullivan Chord diagrams along with their dual representation within
$\Arc$ that appear in this context and in the literature
\cite{CS,C1,C2,CG,TZ}, see \S \ref{chordsection}. When trying to
upgrade them to PROPs one has to allow the gluing of all ``in''
boundaries to all ``out'' boundaries. This poses too many conditions
to make the gluing associative on the topological level.
 As in \cite{cact} this situation can be remedied by weakening the
notion of a PROP to that of a quasi--PROP that is a PROP in which
the associativity holds up to homotopy. The largest
sub--di--operad of $\Diarc$ for which this strategy works is
 $\iooarc$.  For suitable chains the induced chain level structure
is a PROP in the strict sense. The subspace of $\Diarc$ given by
those elements whose  arcs  only run from ``in'' to ``out'' but also
all ``in'' boundaries are hit is called $\Diioarc$.

There is a even a finer structure than the ``in/out'' distinction
which is given by an angle marking. Using angle markings we obtain
an all encompassing rational operad. Here the angles refer to the
angles of the arc graphs which define elements of the $\Arc$ operad
and the marking takes values in $\Z/2\Z$. The space of the elements
of $\Arc$ together with an angle marking is called $\Anarc$. Keeping
with the theme of angle markings we define the space $\Ana$ to be
the space of all elements of $\A$ with additional angle markings.

Armed with these notions we can state a first omnibus theorem summing up essential
results for the topological objects of interest.

\begin{introthm} We have the following topological structures:
\begin{enumerate}

\item The subspace $\Arcno$ of $\Arc$ is a cyclic rational topological  operad. It induces the structure of a cyclic rational topological operad on
the spaces $\Mngn$.

\item $\iooarc$ is a topological quasi--PROP containing $\ioarc$
as a topological sub--quasi--PROP. It induces quasi--PROP structures
on the respective versions of metric Sullivan--Chord diagrams.

\item $\Anarc$ is a cyclic  rational topological operad.
 It ``contains'' all
of the above in the sense that all the structures can be derived
from the gluing in $\Anarc$.

\item $\Diioarc$ is  a quasi--PROP structure extending that of $\ioarc$.
It induces a quasi--PROP structure on the extension of
Sullivan--Chord diagrams for which we define the action on the
Hochschild co--chains in the second part \cite{hoch2}.

\item There exists a CW complex $\Diioarci$ which is homotopy
equivalent to $\Diioarc$. It can be endowed with the structure of
a quasi--PROP that is homotopic to the quasi--PROP structure of
$\Diioarc$.

\end{enumerate}

The same results hold for the restrictions to surfaces with no punctures
and to genus zero surfaces as well as the intersection of these conditions.
\end{introthm}
The main spaces of interest in the sequel will be $\Arcno$ and
$\Diioarc$. As stated above the latter has a CW-model $\Diioarci$.
This model
 is the generalization of $\Cacti^1$ of \cite{cact} to the PROP setting.
 The space $\Arcno$  yields the moduli--space
operations expected for instance from a $D$--brane point of view and
$\Diioarc$ yields String--Topology--type operations, i.e.\ an PROPic
operation of an extension of Sullivan Chord diagrams on the loop
space of a simply connected manifold. Here we use the term Sullivan
Chord diagram for the contracted version Sullivan Chord diagrams;
that is that type of diagrams which one obtains after contracting
ghost edges.

We would also like to mention that there is a generalization
$\Diiooarc$ of $\Diioarc$ where one is also allowed to have arcs
running from ``out'' to ``out''. This is again a quasi--PROP whose
cells give a PROP. It is the generalization of $\iooarc$. It turns
out, see \cite{hoch2}
 that although there is a natural action on the Hochschild
 co-chains,
the  action of this object is not $dg$. There is also a
sub--quasi--PROP $\LDiioarc\subset \Diioarc$ which is the
generalization of the spineless cacti and a sub--CW complex
$\LDiioarci\subset \Diioarci$ corresponding to it. In \cite{hoch2}
we show that there is a natural action of $\LDiioarci$ that
generalizes the $\Xi_2$ action of \cite{MScosimp}.

There is a global object giving rise to all of the needed
structures. This is the CW-complex $\Ana$ of angle marked metric
arc graphs. This space has the structure of a
rational--quasi--modular--operad. Although this is a weak
structure when restricted to the different subspaces which we pick
out, it induces the structures of cyclic operads, rational
operads, quasi-PROPs, etc. we discussed above. In particular this
weak structure is more rigid on the cell level and the correlators
that induce all the different actions are defined for the cells of
$\Ana$ \cite{hoch2}.

The next step is to form cell models for the spaces in question
for which topological structure induce
the honest (not ``quasi'' or ``rational'')
 structures on the chain level.
An essential tool we use is the filtration on $\A$ by the number of
edges of the graph, which gives a filtration of the cells of $\A$.
The gluing rules in $\Arc$ respect this filtration, so it is
possible to consider the associated graded on the cell level. This
for instance gets rid of the codimension one parts which are
responsible for the fact that $\Arcn$ is only  a rational operad.

The essential results used in the sequel for the cell level can be
summed up as follows.

\begin{introthm} We have the following chain level structures:
\begin{enumerate}

\item The relative chains of $\Arc$ in $\A$ form a cellular operad
$\OCArc$ which is filtered by the dimension of the cells.
The associated graded $\GrOCArc$ is also a cyclic operad.

 \item The suboperad $\Gr\OC(\Arcn)$ is a $dg$--operad.
This operad is an operad structure on the collection of graph-complex chain models
of the moduli spaces $\Mngn$ which calculates their cohomology.

\item $CC_*(\A^{\angle})$ is an operad, which contains the relative
chains of $\Gr\OCAnarc$. The chains  $\Gr\OCArc$ are also naturally embedded.

\item $\OC(\Diioarc)$ is a di--operad and a PROP. It is filtered by dimension
and the associated graded $\Gr\OC(\Diioarc)$ is also a di-operad and a PROP.
This is a PROP of the version of Sullivan--Chord diagrams relevant
for our purposes.

\item The cellular chains $CC_*(\Diioarci)$ form a di--operad and
a  PROP --- that is a chain model for  $\Diioarc$. This is a cellular
realization of the above PROP of Sullivan--Chord diagrams.

\end{enumerate}

Again, the same results hold for the restrictions to surfaces with no punctures
and to genus zero surfaces as well as the intersection of these conditions.

\end{introthm}

These cell models are related via the dual graph construction to the
free Abelian groups generated by certain types of ribbon graphs, for
instance: $\Gr\OC(\Arcno)\cong \Rib$, the space of marked ribbon
graphs and $\OC(\ioarc_{\#}) \cong \ppdigraph$ that is the perfectly
partitioned di-graphs which can be thought of as stabilized Sullivan
chord diagrams, see \S\ref{chordsection}.

The paper is organized as follows:

In \S 1, we define the types of graphs and the operations on these
graphs which we will need in the sequel. In \S2 we start be recalling
the constructions of \cite{KLP} albeit in slightly different language
using mainly graphs. We go on to generalize these constructions and
augment the setting by including a filtration on this space. One of the
 new results of this section is the construction of a filtered
cell level operad built on $\Arc$ and its associated graded. The
third paragraph \S3 then details the identification of the
sub-space $\Arcno$ with the moduli space of surfaces with marked
points and tangent vectors at the marked points. The tool here is
the dual graph construction which turns an element of $\Arcno$
into a marked metric ribbon graph. The first main result of the
paper, namely that the associated graded of the cells of $\Arcno$
forms an operad is included in \S4. This paragraph also contains
the identification of this complex with the relevant graph
complex. The last paragraph \S5 contains the generalizations to
the di--operad and PROP structures discussed in the Introduction.
The main link between all the objects are angle marked arc graphs.
The second main result which will be used in the sequel is the
construction of cell and CW models for the graphs corresponding to
Sullivan Chord diagrams whose action will render the String
Topology type operations.

\section*{Acknowledgments} We would like to thank the Max--Planck--Institute
for Mathematics where this work was started, a good portion of
it  was written in the summer of 2005 and the finishing
touches were put on in the summer of 2006. The two papers received their
final form at the MSRI, which we would like to thank for its
hospitality in May 2006. It is a pleasure to thank
 Bob Penner, Ralph Cohen, Jim McClure, Dev Sinha and Craig Westerland for discussions
on various details during various stages of this project.

\section*{Conventions}
We fix $k$ to be a field of arbitrary characteristic. We let $\bar n$ be the set $\set{0}{n}$.

\section{Graphs and Ribbon graphs} \label{Graphs}

\subsection{Classes of Graphs}
In this section, we formally introduce the graphs and the
operations on graphs which we will use in our analysis.

\subsubsection{Graphs} A graph $\Gamma$ is a tuple
$(V_{\Gamma},F_{\Gamma}, \imath_{\Gamma}: F_{\Gamma}\rightarrow
F_{\Gamma},\del_{\Gamma}:F_{\Gamma} \rightarrow V_{\Gamma})$ where
$\imath_{\Gamma}$ is an involution $\imath_{\Gamma}^2=id$ without
fixed points. We call $V_{\Gamma}$ the vertices of $\Gamma$ and
$F_{\Gamma}$ the flags of $\Gamma$. The edges $E_{\Gamma}$ of
$\Gamma$ are the orbits of the flags under the involution
$\imath_{\Gamma}$. A directed edge is an edge together with an
order of the two flags which define it. In case there is no risk
of confusion, we will drop the subscripts $\Gamma$. Notice that
$f\mapsto (f,\imath(f))$ gives a bijection between flags and
directed edges.

We also call $F_{\Gamma}(v):=\del^{-1}(v)\subset F_{\Gamma}$ the
set of flags of the vertex $v$ and call $|F_v({\Gamma})|$ the
valence of $v$ and denote it by $\val(v)$. We also let
$E(v)=\{\{f,\imath(f)\}|f\in F(v)\}$ and call these edges the
edges incident to $v$.

The geometric realization of a graph is given by considering  each
flag as a half-edge and gluing the half-edges together using the
involution $\imath$. This yields a one-dimensional CW complex
whose realization we call the realization of the graph.

\subsubsection{Ribbon graphs}

A ribbon graph with tails is a connected graph together with a
cyclic order of the set of flags of the vertex $F_{\G}(v)$ for
every vertex $v$. A ribbon graph with tails that satisfies
$\val(v)\geq 2$ for all vertices $v$ will simply be called a
ribbon graph.   Notice that we do {\em not} fix $\val(v)\geq3$. We
will call a ribbon graph stable if it does satisfy this condition.

For a ribbon graph with tails, the tail vertices are
$V_{tail}=\{v\in V_{\G}|\val(v)=1\}$, the tail edges
$E_{tail}(\G)$ are the edges incident to the tail vertices  and
the tail flags $F_{tail}(\G)$ are those flags of the tail edges
which are {\em not} incident to the tail vertices.

A  graph with a cyclic order of the flags at each vertex gives
rise to bijections $\Cyc_v:F_v\rightarrow F_v$ where $\Cyc_v(f)$
is the next flag in the cyclic order. Since $F=\amalg F_v$ one
obtains a map $\Cyc:F\rightarrow F$. The orbits of the map
$N:=\Cyc \circ \imath$ are called the cycles or the boundaries of
the graph. These sets have the induced cyclic order.

Notice that each boundary can be seen as a cyclic sequence of
directed edges. The directions are as follows. Start with any flag
$f$ in the orbit. In the geometric realization go along this
half-edge starting from the vertex $\del(f)$, continue along the
second half-edge $\imath(f)$ until you reach the vertex
$\del(\imath(f))$ then continue starting along the flag
$\Cyc(\imath(f))$ and repeat.

An angle is a pair of flags $(f,\Cyc(f))$, we denote the set of
angles by $\angle_{\G}$. It is clear that $f\mapsto (f,\Cyc(f))$
yields a bijection between $F_{\G}$ and $\angle_{\G}$. It is however
convenient to keep both notions.

\subsubsection{The genus of a ribbon graph and its surface} The
genus $g(\Gamma)$ of a ribbon graph $\Gamma$ is given by
$2-2g(\Gamma)=|V_\Gamma|-|E_{\Gamma}|+\C(\G)=\chi(\G)-\C(\G)$
where $\C(\G)=\#cycles$.

The surface $\Sigma(\Gamma)$ of a ribbon graph $\Gamma$ is the
surface obtained from the realization of $\Gamma$ by thickening
the edges to ribbons. I.e.\ replace each 0-simplex $v$ by a closed
oriented disc $D(v)$ and each 1-simplex $e$ by $e\times I$
oriented in the standard fashion. Now glue the boundaries of
$e\times I$ to the appropriate discs in their cyclic order
according to the orientations. This is a surface whose boundary
components are given by the cycles of $\G$. The graph $\Gamma$ is
naturally embedded as the spine of this surface $\Gamma\subset
\Sigma(\Gamma)$. Let $\bar \Sigma(\G)$ be the surface obtained
from $\Sigma(\G)$ by filling in the boundaries with discs. Notice
that the genus of the $\bar \Sigma(\Gamma)$ is $g(\Gamma)$ and
$\chi(\G)=2-2g(\Sigma(\G))$.

\subsubsection{Marked ribbon graphs}
\begin{df}
A ribbon graph together with a distinguished cycle $c_0$ is called
{\em treelike} if
\begin{itemize}
\item[i)] the graph is of genus $0$  and \item[ii)] for all flags either
$f\in c_0$ or $\imath(f)\in c_0$ (and not both).
\end{itemize}
In other words each edge is traversed exactly once by the cycle $c_0$.
Therefore there is a cyclic order on all (non--directed) edges,
namely the cyclic order of $c_0$.
\end{df}

\begin{df}
 A {\em marked ribbon graph} is a
ribbon graph together with a map $\mk:\{cycles\} \rightarrow
F_{\Gamma}$ satisfying the conditions
\begin{itemize}
\item[i)] For every cycle $c$ the directed edge $\mk(c)$ belongs
to the cycle.

\item[ii)] All vertices of valence two are in the image of $\mk$,
that is $\forall v,\val(v)=2$ implies  $v\in Im(\del\circ\mk)$.
\end{itemize}
\end{df}

 Notice that on a marked treelike ribbon graph there is a linear
order on each of the cycles $c_i$. This order is defined by
upgrading the cyclic order to the linear order $\prec_i$ in which
$\mk(c_i)$ is the smallest element.

\subsubsection{Labelling and marking graphs} A labelling of the edges of a
graph $\Gamma$ by a set $S$, is a map
$E_{\Gamma}\rightarrow S$. A labelling of a ribbon graph $\Gamma$
by a set $S$ is a map $\lab\{$cycles of $\Gamma\}\rightarrow S$,
we will write $c_i:=\lab^{-1}(i)$. An angle marking by $S$ is a
map $\amark:\angle_{\G}\rightarrow S$.

\begin{nota}
\label{ribnota} We will write $\Rib_{n,g}$ for the set of marked
ribbon graphs of genus $g$ with $n$ boundary cycles and, by abuse
of notation, also for the free Abelian group generated by this
set.

We set $\Rib:=\amalg_{n,g} \Rib_{n,g}$, and we will again not
distinguish in notation between the set $\Rib$, the free Abelian
group generated by it, and the set $\{\amalg_g \Rib_{n,g}:n\in
\mathbb{N}\}$ to avoid unnecessary clutter.  We also write
$\Rib(n)$ for the set of marked ribbon graphs with $n+1$ cycles
together with a labelling by $\set{0}{n}$ of these cycles. Again
we also denote the free Abelian group generated by this set as
$\Rib(n)$. Finally, to streamline the notation, we will denote the
collection $\{\Rib(n)|n\in \mathbb{N}\}$ simply by $\Rib$.

The meaning of the symbols will always be clear from the context.
\end{nota}

\subsubsection{Spineless marked ribbon graphs}
\label{spinlessgraph} A marked treelike ribbon graph is called
{\em spineless}, if

\begin{itemize}
\item[i)] There is at most one vertex of valence $2$. If there is
such a vertex $v_0$ then $\del(mk(c_0))=v_{0}$.

\item[ii)] The induced linear orders on the $c_i$ are (anti--)compatible
with that of $c_0$, i.e.\ $f\prec_i f'$ if and only if
$\imath(f')\prec_0 \imath(f)$.
\end{itemize}

\subsection{Operations on graphs}
In this section, we will give the basic definitions
of the operations on graphs which we will need.

\subsubsection{Contracting edges} The contraction $\G/e=(\bar
V_{\Gamma}, \bar F_{\Gamma},\bar \imath,\bar \del)$ of a graph
$\G=(V_{\Gamma},F_{\Gamma},\imath,\del)$ with respect to an edge
$e=\{f,\imath(f)\}$ is defined as follows. Let $\sim$ be the
equivalence relation induced by $\del(f)\sim\del(\imath(f))$. Then
let $\bar V_{\Gamma}:=V_{\Gamma}/\sim$, $\bar
F_{\Gamma}=F_{\Gamma}\setminus\{f,\imath(f)\}$ and $\bar \imath:
\bar F_{\Gamma}\rightarrow \bar F_{\Gamma}, \bar\del: \bar
F_{\Gamma}\rightarrow \bar V_{\Gamma}$  be the induced maps.

For a ribbon graph the cyclic order is the one which descends
naturally.

 For a marked  ribbon graph,
we define the marking of $(\bar V_{\Gamma}, \bar F_{\Gamma},\bar
\imath,\bar \del)$ to be $\overline{\mk}(\bar
c)=\overline{\mk(c)}$ if $\mk(c)\notin\{f,\imath(f)\}$ and
$\overline{\mk}(\bar c)=\overline{N\circ \imath(\mk (c))}$ if
$\mk(c)\in\{f,\imath(f)\}$, viz.\ the image of the next flag in
the cycle.

If there is an angle marking, set $f'=N^{-1}(f), f''=\Cyc(f),
g'=N^{-1}(\imath(f))$ and $g''=\Cyc(\imath(f))$, let
$\amark(f',f)=a,\amark(f,f'')=b, \amark(g',\imath(f))=c$ and
$\amark(\imath(f),g'')=d$, after the contraction we set
$\amark(f',g'')=\overline {\bar a\bar d}$ and
$\amark(g',f'')=\overline {\bar b \bar c}$, where we use the
notation $\bar a=1-a \in \Zz$.

\subsubsection{Deleting edges} \label{deledge} The graph
$\G\setminus e=(\bar V_{\Gamma}, \bar F_{\Gamma},\bar \imath,\bar
\del)$ obtained by deletion of an edge $e$  of a graph
$\G=(V_{\Gamma},F_{\Gamma},\imath,\del)$ with respect to an edge
$e=\{f,\imath(f)\}$ is defined as follows:
Set $\bar V_{\G}=V_{\G}, \bar F_{\G}=F_{\G}\setminus
\{f,\imath(f)\}, \bar \imath=\imath|_{\bar F_{\G}},\bar \del
=\del|_{\bar F_{\G}}$. That is delete the edge.
Notice there might be left--over
lone vertices if $e$ was the only edge of the respective vertices,
also the graph might become disconnected.

If the graph was a ribbon graph
the resulting graph again is a ribbon graph using the cyclic order that descends
naturally.  In the case of a marked ribbon graph,
if there was a marking and say $f=\mk(c_i)$ for $f\in e$ then we
set $\mk(c_i)=N^{-1}(f)$ to be the previous flag, otherwise, the marking also
descends naturally.

If there is an angle marking and $f\in e$ with $\val(\del(f))>1$,
 with markings $\amark(f',f)=a,
\amark (f,f'')=b$ for $f'=N^{-1}(f)$ and $f''=\Cyc(f)$ then in
$\G\setminus e$ we set $\amark(f',f'')=\overline {\bar a\bar b}$.

\subsection{Spaces of graphs with metrics}
\subsubsection{Graphs with a metric}

A metric $w_{\Gamma}$ for a graph is a map $E_{\Gamma}\rightarrow
\mathbb{R}_{>0}$. The (global) re-scaling of a metric $w$ by
$\lambda$ is the metric $ \lambda w: \lambda w(e)=\lambda w(e)$.
The length of a cycle $c$ is the sum of the lengths of its edges
$length(c)=\sum_{f\in c} w(\{f,\imath(f)\})$. A metric for a
treelike ribbon graph is called normalized if the length of each
non--distinguished cycle is $1$.
We will write $\MRib_{n,g}$ for the set of metric
marked ribbon graphs of genus $g$ with $n$ boundary cycles.

\subsubsection{Projective metrics}
Notice that there is an $\mathbb{R}_{>0}$--action on $\MRib$ which
scales the metric $\mu$ by an overall factor. This action of
course preserves the genus and number of boundaries. We set
$\PRib:=\MRib/\mathbb{R}>0$ using the same conventions as in
Notation \ref{ribnota}. The elements of $\PRib$ are called graphs
with a projective metric. Notice that one can always choose a
normalized representative for any projective metric.
We set
$\PRib_{n,g}=\MRib_{n,g}/\mathbb {R}_{>0}$.

\subsubsection{The space of metric ribbon graphs}
\label{mribnota}  We endow these above sets
with a  topology by constructing $\PRib(n,g)$  in the standard fashion.
That is we realize them as a subspace of the
 quotient of the disjoint union
of simplices by an equivalence relation.  For each graph $\G\in
\Rib(n,g)$ with $|E(\G)|=k+1$ we fix a $k$-simplex $\Delta_{\G}$.
Using barycentric coordinates for this simplex, a point of this
simplex can be identified with a choice of projective weights on
the edges.
 The points of $\PRib_{n,g}$ can thus be identified with the interior of the disjoint
union over all $\Delta_{\G}:\G\in \Rib_{n,g}$. Furthermore the faces of
$\Delta_{\G}$ correspond to the edges of $\G$. Now,
we use the following identifications: A face of $\Delta_{\G}$ is identified
with $\Delta_{\G/e}$ if $\G/e\in \Rib_{n,g}$.
We give the resulting space the quotient topology (this is actually a CW complex) and identify $\PRib$ with the image of the interiors of the $\Delta_{\G}$.
 Then we
give $\MRib:=\PRib\times\mathbb{R}_{>0}$ the product topology.

\subsubsection{Marked ribbon graphs  with metric and maps of
circles.} For a  marked ribbon graph with a metric, let $c_i$ be
its cycles, let $|c_i|$ be their image in the realization and let
$r_i$ be the length of $c_i$. Then there are natural maps
$\phi_i:S^1\rightarrow |c_i|$ which map $S^1$ onto the cycle by
starting at the vertex $v_i:=\del(\mk(c_i))$ and going around the
cycle mapping each point $\theta\in S^1$ to the point at distance
$\frac{\theta}{2\pi}r_i$ from $v_i$ along the cycle $c_i$. This
observation connects the current constructions to those involving
a more geometric definition of $\Cacti$ in terms of configurations
of circles \cite{Vor,cact} and other geometric constructions
involving such configurations such as the map $\Loop$ used for the
$\Arc$ operad \cite{KLP}. In particular the treelike ribbon graphs correspond
to $\Cacti$ and the spineless treelike ribbon graphs correspond to $\Cact$.

\subsection{Di-graphs and Sullivan Chord diagrams}
\label{chordsection}

\subsubsection{Ribbon Di-graphs} A ribbon graph a di-graph is a
ribbon graph $\G$ together with a $\Zz$ labelling of the cycles of
$\G$: $\io:\{$cycles of $\G\}\rightarrow \Zz$. We call the cycles
$\io^{-1}(0)=:\Out_{\Gamma}$ the outgoing ones and
$\io^{-1}(1)=:\In_{\G}$ the incoming ones. A di-graph is said to
be of type $(n,m)$ if $|\In_{\G}|=n$ and $|\Out_{\G}|=m$. We will
denote the set of these graphs by $\Digraph$.

A ribbon di-graph is called perfectly partitioned, if
$\io(\imath(f))=1-\io(f)$ for every flag  $f$. That is each edge
is part of one input and one output cycle. We will call the set of
these graphs $\ppdigraph$.

A $(S_1,S_2$)-labelled ribbon di-graph is  a ribbon di-graph
together with bijective maps $:\In\rightarrow S_1$ and
$\Out\rightarrow S_2$. We denote the induced map on
$\In\amalg\Out$ by  $\lab$. If $(S_1,S_2)$ is not mentioned, we
will use $S_1=\overline{n}$ and $S_2=\overline{m}$ as the default
indexing sets for a graph of type $(n,m)$.

\subsubsection{Sullivan Chord  and Ribbon Diagrams} There are many
definitions of Sullivan chord diagrams in the literature
\cite{CS,C1,C2,CG,TZ}. We will use the following conventions.

\begin{df}
A  Sullivan chord diagram  is a marked labelled ribbon di-graph
which satisfies the following condition:

\begin{itemize}
\item[i)] after deleting the edges of the incoming cycles one is
left with a forest, i.e.\ a possibly disconnected set of
contractible graphs.

\end{itemize}
We denote the set of these graphs by $\Sulchord$ and call them
Sullivan Chord diagrams.

A {\em strict} Sullivan Chord diagram also satisfies the condition
\begin{itemize}
\item[ii)] It is possible to disjointly embed the incoming cycles
as circles into the plane.
\end{itemize}
We denote this set of graphs by $\StSulchord$.
\end{df}

\begin{rmk}
Notice that in a Sullivan chord diagram there are two types of
edges, those which belong to the $\In$ cycles and those which
belong to the tree part. These are also sometimes called ghost
edges. The edges of the $\In$ cycles are traversed exactly once by
the $\In$ cycles and once by the $\Out$ cycles and the edges.
That is if a flag $f$ is an element of an $\In$ cycle, then
$\imath(f)$ is an element of an $\Out$ cycle. The ghost edges are
traversed in either direction by $\Out$ cycles, that is both $f$
and $\imath(f)$ are elements of $\Out$ cycles.

 If one contracts the ghost edges, one obtains a perfectly
partitioned di-graph. For this reason perfectly partitioned
di-graphs are sometimes called reduced Sullivan chord diagrams.
Notice though that such a diagram is not a strict Sullivan chord
diagram in the above definition, since the $\In$-cycles will share
vertices.
\end{rmk}

\subsubsection{Spaces of di-graphs} To each of these classes of graphs
there is the corresponding space of graphs whose elements are the
graphs of the given type together with a metric. The notation for
these spaces of graphs which we will use is to write $\M$ in front
of the symbol of the graphs, e.g.\ $\M\StSulchord$ and
$\M\ppdigraph$. These spaces have a natural topology. First to each
set of discrete data, such as the labelling and the di-graph
labelling we associate a component, and then forgetting this extra
data each of these components can be identified with a subspace of
$\MRib$ defined by the underlying metric graphs.

 It is clear that the spaces corresponding to $\ppdigraph$ and $\Sulchord$ are
homotopy equivalent. The homotopy is given by changing the metric,
by  homogeneously scaling all the lengths of the ghost edges to
zero.

We will not use strict Sullivan chord diagrams in this paper and
we included them only to make contact with the literature where
they do sometimes appear. The following important fact which one
can show with some effort characterizes the weak homotopy type of
strict Sullivan Chord diagrams.

\begin{claim}
\label{weaksul} The map $p:\M\StSulchord \rightarrow \M\ppdigraph$
which contracts all ghost edges is a weak fibration with
contractible fibers. Hence $\M\ppdigraph$ and $\M\StSulchord$ are
weakly homotopy equivalent.
\end{claim}

\section{The $\Arc$ operad}

\label{arcsection} In this section, we start by giving a brief
review of the salient features of the $\Arc$ operad of \cite{KLP}
which is reasonably self-contained. The presentation of the material
closely follows the Appendix B of \cite{cact}. For full details, we
refer to \cite{KLP}. In addition this review, we furthermore
introduce an equivalent combinatorial language which will be key for
the following, in particular for \cite{hoch2}. Simultaneously, we
introduce new cell level structures and the go on to define new cell
level operads and extensions of the $\Arc$ operad structure.

\subsection{Spaces of graphs on surfaces}
Fix an oriented surface $F_{g,r}^s$ of genus $g$ with $s$
punctures and $r$ boundary components which are labelled from $0$
to $r-1$, together with marked points on the boundary, one for
each boundary component. We call this data $F$ for short if no
confusion can arise.

The piece of the $\Arc$ operad supported on $F$ will be an open
subspace of a space $\Agrs$. The latter space is a CW complex
whose cells are indexed by graphs on the surface $F_{g,r}^s$ up to
the action of the pure mapping class group $PMC$ which is the
group of orientation preserving homeomorphisms of $F_{g,r}^s$
modulo homotopies that pointwise fix  the set which is the union
of the set of the marked points on the boundary and the set of
punctures. A quick review in terms of graphs is as follows.

\subsubsection{Embedded Graphs}
By an embedding of a graph $\Gamma$ into a surface $F$, we mean an
embedding  $i:|\Gamma|\rightarrow F$ with the conditions

 \begin{itemize}
\item[i)] $\Gamma$ has at least one edge.

\item[ii)] The vertices map bijectively to the marked points on
the boundaries.

\item[iii)] No images of two edges are homotopic to each other, by homotopies
fixing the endpoints.

\item [iv)] No image of an edge is homotopic to a part of the
boundary, again by homotopies fixing the endpoints.
\end{itemize}

Two embeddings are equivalent if there is a homotopy of embeddings
of the above type  from one to the other. Note that such a
homotopy is necessarily constant on the vertices.

The images of the edges are called arcs. And the set of connected
components of $F\setminus i(\Gamma)$ are called complementary
regions.

Changing representatives in a class  yields natural bijections of
the sets of arcs and connected components of $F\setminus
i(\Gamma)$ corresponding to the different representatives.
 We can therefore associate to each equivalence class of embeddings
its sets of arcs together
with their incidence conditions and connected
components --- strictly speaking of course the equivalence classes of these
objects.

\begin{df}
By a graph $\g$ on a surface we mean a triple $(F,\Gamma,[i])$ where
$[i]$ is an equivalence class of embeddings of $\Gamma$ into that
surface. We will denote the isomorphism class of complementary
regions by $\comp(\g)$. We will also set $|\g|=|E_{\G}|$. Fixing the
surface $F$, we will call the set of graphs on a surface
$\graphs(F)$.
\end{df}

\subsubsection{A linear order on arcs}
Notice that due to the orientation of the surface the graph
inherits an induced linear order of all the flags at every vertex
$F(v)$ from the embedding. Furthermore there is even a linear
order on all flags by enumerating the flags first according to the
boundary components on which their vertex lies and then according
to the linear order at that vertex. This induces a linear order on
all edges by enumerating the edges by the first appearance of a
flag of that edge.

\subsubsection{The poset structure}
The set of such graphs on a fixed surface $F$ is a poset. The
partial order is given by calling $(F,\Gamma',[i'])\prec
(F,\Gamma,[i])$ if $\Gamma'$ is a subgraph of $\Gamma$ with the
same vertices and $[i']$ is the restriction of $[i]$ to $\Gamma'$.
In other words, the first graph is obtained from the second by
deleting some arcs.

 We
associate a simplex $\Delta(F,\Gamma,[i])$ to each such graph.
$\Delta$ is the simplex whose vertices are given by the set of
arcs/edges enumerated in their linear order. The face maps are
then given by deleting the respective arcs. This allows us to
construct a CW complex out of this poset.

\begin{df}
Fix $F=F_{g,n}^s$. The space $\A_{g,n}^{\prime s}$ is the space
obtained by gluing the simplices $\Delta(F,\Gamma',[i'])$ for all
graphs on the surface according to the face maps.
\end{df}

The pure mapping class group naturally acts on $\A_{g,n}^{\prime s}$
and has finite isotropy \cite{KLP}.

\begin{df}
The space $\A_{g,r}^s:= \A_{g,r}^{\prime s}/PMC$.
\end{df}

\subsubsection{CW structure of $\A_{g,r}^s$}

\begin{df}
Given a graph on a surface, we call its $PMC$ orbit its arc graph.
If  $\g$ is a graph on a surface, we denote by $\bar{\g}$ its arc
graph or $PMC$ orbit. We denote the set of all arc graphs of a
fixed surface $F$ by $\arcgraphs(F)$. A graph is called exhaustive
if there are no vertices $v$ with $\val(v)=0$. This condition is
invariant under $PMC$ and hence we can speak about exhaustive arc
graphs.  The set of all exhaustive arc graphs on $F$ is denoted by
$\carcgraphs(F)$.
\end{df}

Notice that since the incidence conditions  are preserved, we can
set $|\bar{\g}|=|\g|$ where $\g$ is any representative and
likewise  define $\comp(\bar{\g})$. We call an arc graph
 exhaustive if and only if it contains no isolated vertices, that is
 vertices with $\val(v)=0$.

Now by construction  it is clear that $\Agrs$ is realized as a CW
complex which has one cell for each arc graph $\bar{\g}$ of
dimension $|\g|-1$. Moreover the cell for a given class of graphs is
actually a map of a simplex whose vertices correspond to the arcs in
the order discussed above. The attaching maps are given by deleting
edges and identifying the resulting face  with its image. Due to the
action of $PMC$ some of the faces of might become identified by
these maps, so that the image will not necessarily be a simplex. The
open part of the cell will however be an open simplex. Let $C(\ba)$
be the image of the cell and $\dot C(\ba)$ be its interior, then
\begin{equation}
\A_{g,r}^s=\cup_{\ba\in \arcgraphs(F_{g,r}^s)} C(\ba),\quad
\A_{g,r}^s=\amalg_{\ba\in \arcgraphs(F_{g,r}^s)} \dot C(\ba)
\end{equation}
 Let $\Delta^n$ denote the standard
$n$--simplex and $\dot \Delta$ its interior then
 $\dot C(\g)=\mathbb R^{|E_{\G}|}_{>0}/\mathbb{R}_{>0}=\dot
\Delta^{|E_{\G}|-1}=:C(\G)$ which only depends on the underlying graph
$\G$ of $\g$.

 This also  means that the space $\Agrs$ is
filtered by the cells of dimension less than or equal to $k$. We
will use the notation $(\Agrs)^{\leq k}$ for the pieces of
this filtration.

\subsubsection{Open-cell cell complex}
It is clear by construction that the $\Arc$ operad again has a
decomposition into open cells.

\begin{equation}
\Arc_g^s(n)=\amalg_{\g=\carcgraphs(F_{g,n+1}^s)}\dot C(\g)
\end{equation}
again $\dot C(\g)=\mathbb R^{|E_{\G}|}_{>0}/\mathbb{R}_{>0}=\dot
\Delta^{|E_{\G}|-1}:=\dot C(\G)$ only depends on the underlying
graph $\G$ of $\g$.

\label{graded} We will denote the free Abelian group generated by
the $C(\a)$ as above by $\OCArc_{g}^s(n)$. We will write $\OCArc(n)
=\amalg_{g,s}\OCArc_g^s(n)$ and $\OCArc=\amalg_n \OCArc(n)$. We
choose the notation to reflect the fact that we are strictly
speaking not dealing with cellular chains, however see \S
\ref{graphcompsection}.

$\OCArc(n)$ is also graded by the dimension of the cells, we will
write $\OCArc(n)^k$ for the subgroup generated by cells of dimension
$k$ and we will also write $\OCArc(n)^{\leq k}$ for the subgroup of
cells of dimension $\leq k$. It is clear that $\OCArc(n)^{\leq k}$
induces a filtration on $\OCArc(n)$ and that the associated graded
is isomorphic to the direct sum of the $\OCArc(n)^k$
\begin{equation}
\GrOCArc := Gr(\OCArc(n),\leq) \simeq \bigoplus_k \OCArc^k(n)
\end{equation}

The differential $\del$ of $\A_{g,r}^s$ also descends to $\OCArc$
and $\GrOCArc$ by simply omitting the cells which are not in $\Arc$.
Applying the differential twice will kill two arcs,
 each original summand will either be twice treated as
zero or appear with opposite sign as in $\A_{g,r}^s$. Hence the
differential squares to zero.

\subsubsection{Relative cells}
The complex $\OCArc_{g}^s(n)$ and the isomorphic complex
$\GrOCArc_G^s(n)$ can be identified with the complex of relative
cells $CC_*(A,A\setminus\Arc)$ .

\subsubsection{Elements of the $\Agrs$ as projectively weighted graphs}
Using barycentric coordinates for the open part of the cells the
elements of $\Agrs$ are given by specifying an arc graph
together with a map $w$ from the edges of the graph $E_{\Gamma}$ to
$\mathbb{R}_{>0}$ assigning a weight to each edge s.t.\ the sum of
all weights is 1.

\begin{figure}
\epsfxsize = 4in \epsfbox{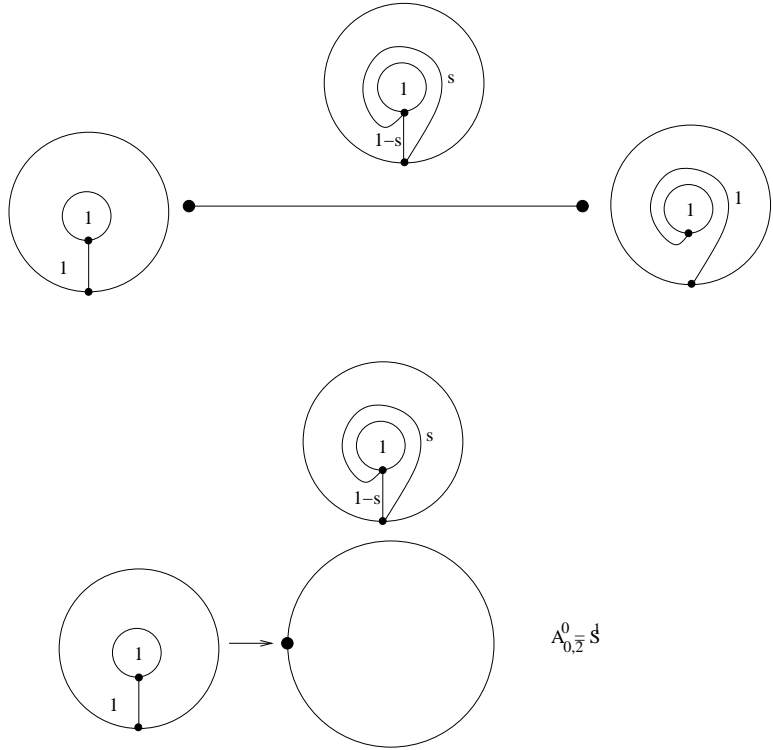} \caption{\label{f02delta}
The space $\A^0_{0,2}$ is given as the CW decomposition of $S^1$
with one $0$--cell and one $1$--cell. It can be thought of as the
quotient of the interval in which the endpoints are identified by
the action of the pure mapping class group. The generator of $CC_*(S^1)$ is
called $\Delta$.}
\end{figure}

Alternatively, we can regard the map $w:E_{\Gamma}\rightarrow
\mathbb{R}_{>0}$ as an equivalence class under the equivalence
relation of, i.e.\ $w\sim w'$ if $\exists \lambda \in
\mathbb{R}_{>0}\forall e\in E_{\Gamma} \;  w(e)=\lambda w'(e)$.
That is $w$ is a projective metric. We call the set of $w(e)$ the
projective weights of the edges.
 In the limit, when the projective weight of an
edge goes to zero, the edge/arc is deleted, see \cite{KLP} for more
details. For an example see Figure \ref{f02delta}, which is
discussed below.

An  element $\alpha\in \Agrs$ can be described by a tuple
$\alpha=(F,\Gamma,\overline{[i]},w)$ where $F$ and $\G$ are as
above, $\overline{[i]}$ is a PMC orbit of an equivalence class of
embeddings and $w$ is a projective metric for $\G$. Alternatively it
can be described by a tuple $(\bg,w)$ where $\bg\in \arcgraphs(F)$
and $w$ is a projective metric for the underlying abstract graph
$\G$.
\begin{ex} $\A_{0,2}^0=S^1$. Up to PMC there is a unique graph
with one edge and a unique graph with two edges. The former gives a
zero--cell and the latter  gives a one--cell whose source is a
1--simplex. Its two subgraphs with one edge that correspond to the
boundary lie in the same orbit of the action of PMC and thus are
identified to yield $S^1$. The fundamental cycle is given by
$\Delta$ of Figure \ref{f02delta}.
\end{ex}

\subsubsection{Drawing pictures for Arcs.} 
\begin{figure}
\epsfxsize = 3in
\epsfbox{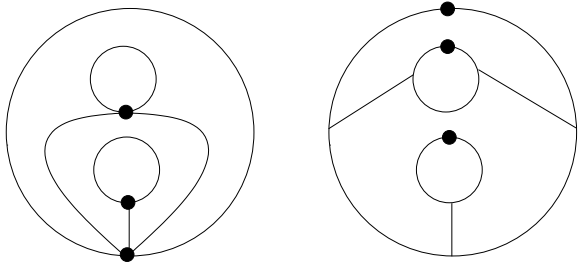}
\caption{\label{different}
An arc graph and its alternate depiction with disjoint arcs not hitting the marked points on the boundary.}
\end{figure}

There are several
pictures one can use to view elements of $\A$. In order
to draw elements them it is useful to expand the
marked point on the boundary to an interval or window, and let the
arcs end on this interval according to the linear order.
Equivalently, one can mark one point of the boundary and let the
arcs end in their linear order anywhere but on this point. We will
mostly depict arc graphs in the latter manner. See Figure
\ref{different} for an example of an arc graph ---all arcs running
to the marked points--- and its alternate depiction with {\em none} of the arcs
hitting the marked point and all arcs having {\em disjoint
endpoints}.

\begin{nota}
Since in the following we will always be dealing with arc graphs,
we will now omit the over-line in the notation. Hence we will
write $\g\in \arcgraphs(F)$. We also fix that $\G(\g)$ is the
underlying graph. Furthermore elements of $\A_{g,r}^s$ will usually
be called $\a$ and $\b$. If $\a\in \A_{g,r}^s$ we fix that
$\g(\a)$, $\G(\a)$ and $w(\a)$ are the underlying arc graph, its
underlying graph and the projective metric, respectively.
\end{nota}

\begin{df}
We define the Euler characteristic of an element $\alpha \in
A_{g,r}^s$ to be $\chi(\alpha)=|\comp(\a)| - |E(\Gamma(\alpha)|$.
\end{df}

\begin{lem}
\label{Euler}
 The following inequality holds
\begin{equation}
\chi(\alpha)\geq \chi(F(\alpha))
\end{equation}
and the equality holds if and only if the complementary regions
are polygons.
\end{lem}

\begin{proof}
If the complementary regions are polygons, we obtain a
triangulation of the surface and hence
$\chi(F)=|V_{\G}|-|E_{\G}|-n+|\comp(\a)|$. There are $|V_{\G}|=n$
vertices and $|E_{\G}|+n$ edges since the boundaries also count as
edges. In this count, the complementary region contributes 1, as
it should for a polygon. If the complementary  regions would have
some topology then their Euler-characteristic would be strictly
less than 1 and the Euler-characteristic of the surface would be
bigger; whence the claim.
\end{proof}

\subsection{Topological operad structure}
\label{gluingpar}
\subsubsection{The spaces $\Arc(n)$}
We begin by reviewing the construction of \cite{KLP}. We then recast it
into a purely combinatorial way. This will allow us to define the actions of
\cite{hoch2} more simply, but also allow us to show that although $\Arcn$ is not
an operad on the topological level, it is a rational operad
and gives rise to a cellular operad.

\begin{df}
We define $\Arc_g^s(n)\subset \A_{g,n+1}^s$ to be the subset of
those weighted arc graphs whose arc graph is exhaustive. We define
$\Arc(n):=\coprod_{s,g\in \mathbb{N}}\Arc_g^s(n)$.\footnote{Unfortunately there is
a typo in the definition of $\Arc(n)$ in \cite{KLP} where
 $\coprod$ was inadvertently replaced by the direct limit.}
\end{df}

Notice that the spaces $\Arc(n)$ have a natural operations of $\Sn$
which permutes the labels $\{1,\dots,n\}$ and one of $\Snn$ which
permutes the labels $\{0,\dots,n\}$. Also notice that the spaces
$\Arc_g^s(n)$ inherit the grading and filtration from $\A_g^s(n)$.
This is also true for their unions $\Arc(n)$ and we will write
$\Arc(n)^{\leq k}$ for these pieces. That is if $\a\in \Arc(n)^{\leq
k}$ then $|E(\G(\a))|\leq k+1$.

\subsubsection{Topological description of the gluing \cite{KLP}}
To give the composite $\alpha \circ_i\alpha'$ for two arc families
$\alpha=(F,\Gamma,\overline{[i]},w)\in \Arc(m)$ and
$\alpha'=(F',\Gamma',\overline{[i']},w') \in \Arc(n)$ one most
conveniently chooses metrics on $F$ and $F'$. The construction does
not depend on the choice. With this metric, one produces a partially
measured foliation in which the arcs are replaced by bands of
parallel leaves (parallel to the original arc) of width given by the
weight of the arc. For this we choose the window representation and
also make the window tight in the sense that there is no space
between the bands and between the end-points of the window and the
bands. Finally, we put in the separatrices. The normalization we
choose is that the sum of the weights at boundary $i$ of $\alpha$
coincides with the sum of the weights at the boundary $0$, we can
also fix them both to be one. Now when gluing the boundaries, we
match up the windows, which have the same width, and then just glue
the foliations. This basically means that we glue two leaves of the
two foliations if they end on the same point.  We then delete the
separatrices.  Afterwards,  we collect together all parallel leaves
into one band. In this procedure, some of the original bands might
be split or ``cut'' by the separatrices. We assign to each band one
arc with weight given by the width of the consolidated band. If arcs
occur, which do not hit the boundaries, then we simply delete these
arcs. We call these arcs or bands ``closed loops'' and say that
``closed loops appear in the gluing''.

\begin{thm}\cite{KLP}
Together with the gluing operations above, the spaces $\Arc$ form
a cyclic operad.
\end{thm}

\begin{lem}
\label{gradinglem} The gluing operations are compatible with the
filtrations, i.e.\ $\Arc(m)^{\leq p}\circ_i \Arc(n)^{\leq q}\subset
\Arc(n+m-1)^{\leq p+q}$.
\end{lem}
\begin{proof}
If $k$ is the number of arcs at the boundary $i$ of $\a$ at the
boundary
 and $l$ is the number of arcs at the boundary $0$ of $\b$ then
the number of glued arcs resulting from the operadic composition
$\circ_i$ is at most $p+q-1$. This comes from the fact that each
``cut'' that is separatrix of the glued foliation contains at least
one of the separatrices of the two foliations before gluing. Hence
the number of arcs of $\a\circ_i\b$
 is at most $p+1-k+q+1-l+k+l-1=p+q+1$ and  the claim follows.
\end{proof}

\subsubsection{Combinatorial description}
\begin{figure}
\epsfxsize = \textwidth \epsfbox{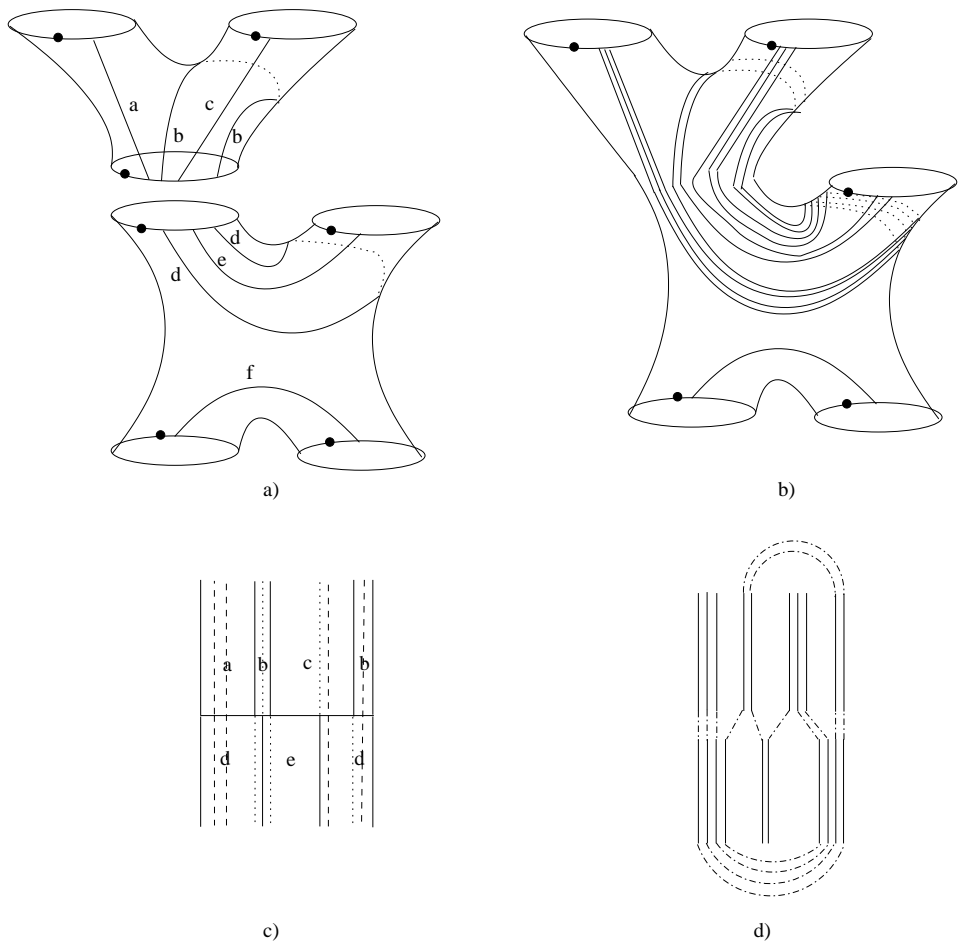}
\caption{\label{gluingex1} a) The arc graphs which are to be glued
assuming the relative weights a,b,c,d and e as indicated by the
solid lines in c). b) The result of the gluing (the weights are
according to c). c) The combinatorics of cutting the bands. The
solid lines are the original boundaries the dotted lines are the
first cuts and the dashed lines represent the recursive cuts. d)
The combinatorics of splitting, and joining flags.}
\end{figure}

One can also give a purely combinatorial description of the
gluing in which we define the new graph on the glued surface.
An example of a gluing is given in Figure \ref{gluingex1}.

In the first step we normalize as above to make
the sums of weights match.  The weights of the arcs
incident to the boundaries $0$ and $i$ then define two partitions $P_0$
and $P_i$ of the unit interval where each flag incident to the
vertices $0$ and $i$ corresponds now to a subinterval of length
given by the weight of the respective arc.

Let's first also assume that there are no arcs running from $0$ to
$0$ or from $i$ to $i$. We will call such arcs recursive. Hence we
are for the moment assuming that there are no recursive arcs.
 Let $\mathcal{P}_{\circ}$ be the
biggest common sub-partition. If in the common partition
$\mathcal{P}_{\circ}$ a subinterval $I'$ of $P_0$ or $P_i$ is
subdivided, let $f$ be the flag associated to $I'$ and duplicate
the edge $e=\{f,\imath(f)\}$ so that there are now two copies
$e^{(1)}$ and $e^{(2)}$ of $e$ which are embedded to be parallel.
Associate the copies $f^{(1)}$ and $f^{(2)}$ of the flag $f$ in
their natural order to the intervals of $P_{\circ}$. Now delete
the vertex $0$ and all flags of the vertex $0$ of $\alpha'$ and
the vertex $i$ and all the flags of the vertex $i$ of $\alpha$.
Then glue the remaining graphs with the duplicated flags by
defining the involution $\imath_{\circ}$ on the flags of this
graph by extending the involution induced by the two original arc
graphs, by setting $u=\imath_{\circ}(d)$ if $\imath(u)$ and
$\imath(d)$ correspond to the same subinterval $I''$ of the
partition $\mathcal{P}_{\circ}$. The weight of such a glued edge
is defined to be the length of the subinterval $I''$.

In the case that there are recursive arcs the combinatorics are a
little more difficult. We start out as above and will call a biggest
common partition $\mathcal{P}$. A subinterval is called recursive if
it belongs to a flag which corresponds to an recursive arc. Now, we
have to complete the partitions. This means that if a subinterval is
recursive, we have to partition the interval corresponding to
$\imath(f)$ in the induced orientation just as the interval
corresponding to $f$. This amounts to cutting the band corresponding
to the edge $(f,\imath(f))$ according to the partition. One might
have to iterate this process. The process will stabilize however,
since there are only finitely many intervals and cuts; and for that
matter only finitely many bands after gluing \cite{KLP}. This will
yield a partition $\mathcal{P}_{\circ}$ of the interval. We now
proceed with the combinatorics as above. Thus by replicating edges,
we obtain two flags per subinterval of $\mathcal{P}_{\circ}$, one on
each side. The glued graph is now defined in three steps. First fix
one vertex for each subinterval an define the two flags incident to
the subinterval to be incident to the vertex. Secondly, delete these
vertices of valence two and their flags. Thirdly glue the remaining
unpaired flags to edges using $\imath_{\circ}$ which is defined
analogously as above.

For the embedding one again uses the window description. After
possible duplications the arcs can be arranged to end in the
mid-point of the subintervals of $P_{\circ}$ and the embedding is
defined by connecting two arcs if their endpoints coincide.

In both cases, there is a last stabilization step, if the edges
corresponding to the first and last interval are parallely
embedded, merge them to one edge and define the weight of this
edge to be the sum of the weights of the parallel edges.

Notice that  this prescription automatically deletes any closed
arcs, viz.\ bands not hitting the boundaries, that might appear in
the second step.

\begin{nota}
\label{splitnota}
 We wish to fix the following terminology.
In the gluing of two families there are two possibilities for a
non--recursive arc which is glued. It is either {\em spilt} if the
flag not incident to the boundary on which the gluing is performed
is duplicated or {\em not split} if the opposite is the case. If
the flag is duplicated, it gives rise to a new angle in the glued
graph which we call a split angle.

We also call an arc family {\em twisted at the boundary $i$} if
the first and the last
 arcs incident to the boundary $i$ become
parallel if one allows homotopies moving the endpoints around the
boundary $i$. We use the same terminology for the underlying
graphs.
\end{nota}

\subsubsection{Combinatorics of the gluing}
\label{gluingcomb} For future reference, we wish to detail the
combinatorics of the gluing if the weights on the arcs on both the
boundaries which are to be glued vary freely.

First we treat the case where there are no arcs running from $0$ to
$0$ or from $i$ to $i$. Also for the moment assume not both
boundaries are twisted.
 Fix $\a\in \Arc(n)$ and $\b\in \Arc(m)$
and then let $\a$ and $\b$ vary freely inside their cell. Assume
that inside the cell $\a$ has $k$ arcs with freely varying weights
on the boundary $i$ and $\b$ has $l$ arcs with freely varying
weights on the boundary $0$. Then the gluing will generically glue
these $k+l$ arcs into $k+l-1$ arcs, since the generic common
partition will have exactly that many components. The combinatorics
are hence the decomposition of $\Delta^{k-1} \times
\Delta^{l-1}=\sum_{\s\in Sh(k,l)} \Delta^{k+l-2}$ where $Sh(k,l)$
are the shuffles of type $k,l$. A non--generic situation happens if
some of the partition points match up this corresponds passing to
the faces of the simplices. This set is at least of codimension one.

If both gluing  boundaries are twisted, we generically only obtain
$k+l-2$ arcs since the two outside arcs will always be parallel
after gluing.

Now say that  there are $s$ arcs running from $0$ to $0$ and $t$
arcs  from $i$ to $i$. Let $k$ and $l$ be as above. We notice that
the number of arcs after gluing is again generically $k+l-1$ arcs.
For this, we notice that the final number of arcs is given by the
half the number of replicated flags for the gluing minus the number
of flags which run from the boundary back to itself. Or in the
notation of \S\ref{gluingcomb}, the number of subintervals of
$\mathcal{P}_{\circ}$ minus half the number of recursive flags. So
we are left with one half the number of the non--recursive flags.
The number of these flags corresponds to the interval markings or
cuts on the respective subintervals. Now counting only these
markings, we count each original cut (subinterval of $\mathcal{P}$
before the iteration process of \S\ref{gluingcomb}) exactly once,
since the iteration stops exactly if the opposite side of the
interval is non--recursive. The number of these cuts is $k+l$ as
above. The indexing of the occurring combinatorics can again be done
by $Sh(k,l)$ as per the original cuts, viz.\ $\mathcal{P}$, but the
topology of the resulting arcs is more complex. Now there are two
types of non--generic situations. The first is as above, that is the
situation corresponding to the partitioning points matching up. This
is again of codimension one. The second is given by bands which form
closed loops that are erased. This is actually even of codimension
two, since  in order for a closed loop to form two of the cuts have
to align.

\subsubsection{Partial operad structure on $\A^s_{g,r}$}
We wish to point out that the gluing we defined in \cite{KLP} that
is reviewed above extends to a partial cyclic operad structure on
$\A(n):=\amalg_{s,g}A^s_{g,n+1}$. The $\Snn$-action is taken to be
the one permuting the labels. The gluings $\circ_i$ above are well
defined as soon as there are arcs incident to {\em both} $i$ and $0$
on the surfaces which are glued along these boundaries.

\subsection{Operad Structures on the Cell level}
\label{cellops}
\subsubsection{Singular Chains/PL Chains}
First it is clear that there are family operations which have as
input the singular co-chains of the $\Arc$ operad and as output
products of these chains. So that if we move to PL-chains, we
indeed get an operad action on the chain level. This is the point
of view taken in \cite{KLP} where we showed that this cell level
operad has a Gerstenhaber-Batalin-Vilkovisky (GBV) structure up to
explicitly given homotopies.

\subsubsection{Operad structures on the free Abelian group of open
cells}
\label{cellop}

 There are two operad structures on the open/relative cell level.
 The first is the more na\"{\i}ve induced operad structure on the open
 cells. And the second one is its associated graded which is more suited for the description
 of moduli space.

To give the first operad structure we claim that as sets
\begin{equation}
\label{setcirc}
\dot C(\a)\circ_k \dot C(\b)=\amalg_{i\in I}\dot C(\g_i)
\end{equation}
for some index set $I$. This reason for this is that if a point of a
cell $C(\g)$ is in the image, the whole cell is. It is possible to
obtain all weights of by varying the weights of the $\a$ and $\b$
accordingly, see \S \ref{gluingcomb} and Lemma \ref{bijectionlem}
below. By Lemma \ref{gradinglem}
 the dimension of the
 cells $C(\g_i)$ appearing in the image
  is less or equal to the sum of the dimensions $\dot C(\a)$ and
  $\dot C(\b)$. Finally, again by \S \ref{gluingcomb}
  and  Lemma \ref{bijectionlem}, we have cells of the
 top-possible dimension if and only if not both the boundaries on which one glues
are twisted. Lastly, the gluing map between cells is also 1-1 precisely
when there is no twisting at both boundaries and there are no closed loops.

\begin{lem}
\label{bijectionlem}
 The
map $\circ_i$ restricted to
$\dot C(\a)\times\dot C(\b)$ is  a bijection onto its
image if and only if not  both $\a$ is twisted
at the boundary $i$ and $\b$ is twisted at the boundary $0$ and there are no closed loops.
In the case not both boundaries are twisted $\circ_i$ restricted
to the set of metric arc graphs for which no closed loops appear is bijective onto is image.
\end{lem}

\begin{proof}
Fix two surfaces $F_{g,r}^s$ and $F_{g',r'}^{\prime \, s'}$ and let
$F=F_{g+g',s+s'}^{\prime \prime \, s+s'}$
 the surface obtained from $F$ and $F'$
by gluing the boundary $0$ of $F'$ to the boundary $i$ of $F$ and
let $l$ the simple closed curve which is the image of the
boundaries $i$ and $0$ in the glued surface. To avoid too many
super- and subscripts set $\Arc(F):=\Arc(F_{g,r}^s):= \A_{g,s}^s$
and likewise for $F'$ and $F''$. Assuming that there
are no closed loops, the gluing procedure is such,
that the preimage of $\circ_i$ from $\Arc(F)\times
\Arc(F')\rightarrow \Arc(F'')$ can be described as follows. Given
a weighted arc graph $\g$ on the glued surface $F''$ all the
preimages that is $\a$, $\b$ such that $\a\circ_i\b=\g$ are
obtained by first fixing a point $p$ on the curve $l$, then
cutting the surface $F''$ along $l$ and finally merging all the
 parallel bands which might occur after cutting while
summing their weights. Here we allow the point $p$ to ``split an
arc''. By this we mean that given a fixed arc with weight $w$ we
draw a parallel arc  to this arc (say to the right of the arc) and
choose $w_1,w_2$ with $w_1+w_2=w$ and  put a point between the two
parallel arcs, that is to the right of the original arc. Fixing
the cells $\dot C(\a)$ and $\dot C(\b)$ we see that the choice of
the point on $l$ is fixed up to moving the point along $l$ but not
crossing any arc. From this description it follows that $\circ_i$
is injective precisely if this point $p$ does not ``split an
arc''. In the case of closed loops, we see that we cannot detect them on the glued
side and that accordingly, the parameters given by their width are free parameters
in the pre--image.
\end{proof}

We will accordingly split the indexing set $I$ of equation
\ref{setcirc} into $I'$ which indexes the cells that are in the
image of the gluings which do not exhibit any closed loops and $I''$
which indexes the cells whose graphs are obtained by erasing closed
loops. As sets, we have:
\begin{equation}
\label{setcirctwo}
\dot C(\a)\circ_k \dot C(\b)=(\coprod_{i\in I'}\dot C(\g_i')) \amalg
(\coprod_{i''\in I''}\dot C(\g_i') )
\end{equation}

We define:
\begin{equation}
\label{setop} \dot C(\a)\tilde\circ_k \dot C(\b):=\begin{cases}\sum_{i'\in
I'}\pm\dot C(\g_{i'})&\begin{tabular}{l} if {\em not} both $\a$ is twisted
at the boundary $k$\\ and $\b$ is twisted at the boundary $0$
\end{tabular}\\
0&\begin{tabular}{l} if  both $\a$ is twisted
at the boundary $k$\\ and $\b$ is twisted at the boundary $0$
\end{tabular}
\end{cases}
\end{equation}
in $\OCArc(n)$ where the sign comes from the orientation of the
cells given by the enumeration of the arcs.

\begin{lem}
With the induced operad structure $\OCArc$ is a cyclic operad which
respects the filtration by dimension.
\end{lem}

\begin{proof}
The fact that the $\Sn$-actions permuting the labels together with
the operations defined in equation (\ref{setop}) yield an operad
structure follows from the observations of the previous paragraph
summed up \S \ref{gluingcomb}. The cyclicity under the $\Snn$
action is inherent. The last statement follows from Lemma
\ref{gradinglem}.
\end{proof}

\begin{cor}
The set of associated graded spaces $\GrOCArc$ together with the
action of the permutation groups on the labels and the operations
induced from equation (\ref{setop}) form a cyclic operad.
\end{cor}
\qed

\begin{figure}
\epsfxsize = \textwidth
\epsfbox{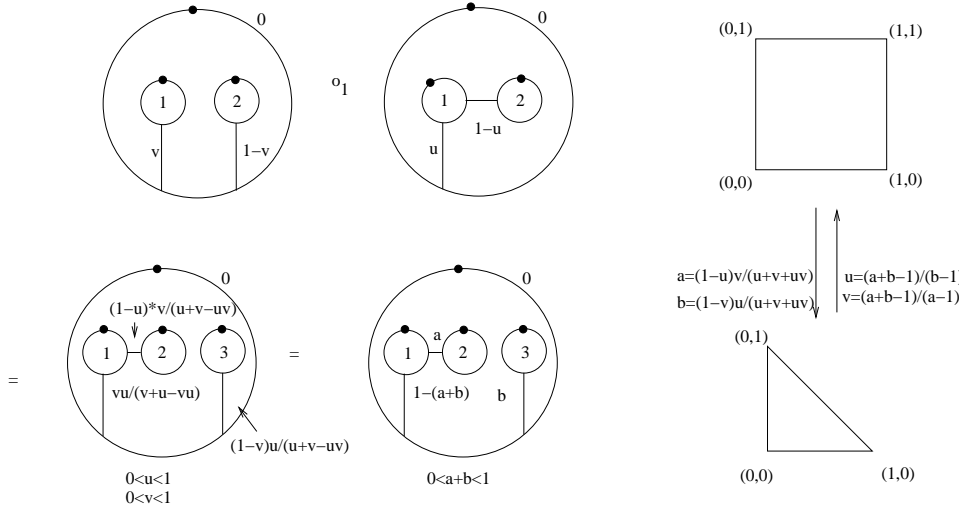}
\caption{\label{boundprobs3}
 A gluing in $\OCArc$ together with a reparameterization of $(0,1)\times(0,1)$
to the inside of the 2--simplex}
\end{figure}

\begin{rmk}
\label{nodg} Both operad structures are not  structures of
$dg$-operads, for the following simple reason. If we glue together
two surfaces then  on the glued side, the limit in which the
weight of all the arcs hitting the separating curve which is the
image of the two glued boundaries goes to zero is possibly allowed
and possibly contributes to the boundary of the cell. This limit,
however, is not allowed for the two components, i.e.\ the result
of this limit does not lie in $\del \dot C(\a) \circ_i \dot \b
\cup \dot C(\a) \circ_i \del \dot \b $. For an example, see Figure
\ref{boundprobs3}. Here the two families have trivial boundary
individually in $\OCArc$, but their composition has a non--trivial
boundary, since the limit where $a$ tends to zero is well defined
in $\Arc$.

We wish to point out that this limit is not allowed in $\Arcn$, see
 section \S \ref{extasec}.
On the cellular level there are ways to remedy this situation,
see \S \ref{subspaces}.

\end{rmk}

Let us make the cell level operation on the associated graded
explicit. If $C(\a)\in \Arc^k(n)$ and $C(\b)\in \Arc^l(m)$  we
obtain
$$
\dot C(\a)\tilde\circ_k \dot C(\b)=\sum_{i\in I}\pm\dot C(\g_i) +
\sum_{j\in J} \pm C(\delta_j)
$$
for some index sets $I$ and $J$ with $C(\g_s)\in \Arc^{k+l}(n+m)$
and $\C(\delta_j)\in \Arc^{<k+l}(n+m)$.

 Therefore, we get an induced operad structure on $\GrOCArc$
 by
setting
\begin{equation}
\label{gradedcellop} \dot C(\a)\circ_k \dot C(\b)=\sum_{i\in I}\pm
C(\g_i)
\end{equation}
where $\pm$ is the usual sign corresponding to the orientation,
which is obtained from the shuffle of the edges that is induced by
the respective shuffle in the product of sub-simplices
$\Delta^{k-1}\times \Delta^{l-1}$ as discussed in
\S\ref{gluingcomb}.

The  result of the operation (\ref{gradedcellop}) is zero exactly
if the index set $I$ is empty and this is the case if and only if
both $\a$ is twisted at $i$ and $\b$ is twisted at $0$. Otherwise
$I$ coincides with the set $Sh(k,l)$ of $(k,l)$-shuffles.

\begin{figure}
\epsfxsize = \textwidth
\epsfbox{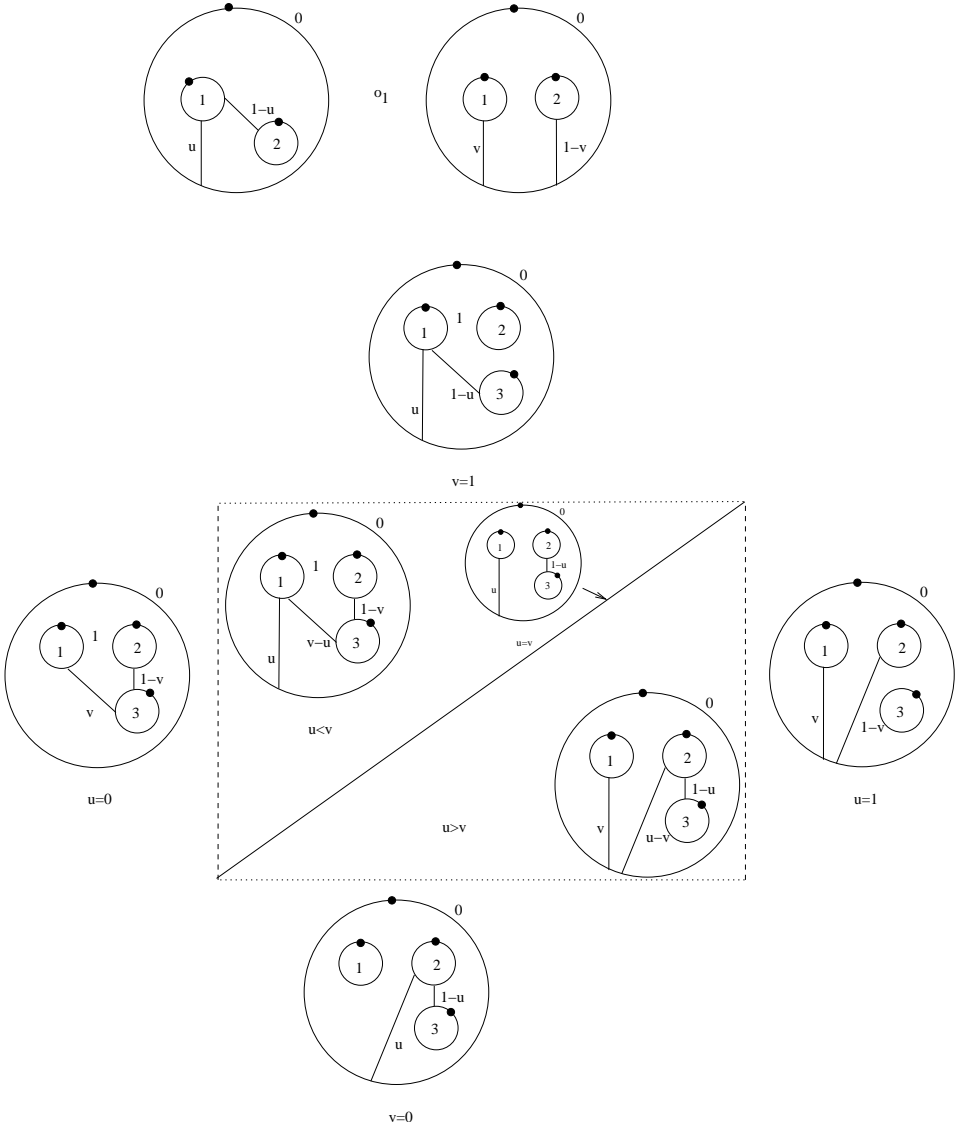}
\caption{\label{boundprobs}
An example of gluing cells. The dashed lines denote the closure of the open cells and the smaller arc
graph denotes the codimension one ``diagonal'' family}
\end{figure}

An example of the gluing is given in Figure \ref{boundprobs}. Here the ``diagonal'' family is of codimension 1 and is included in the open cell gluing while in the graded gluing it is set to zero.

\subsubsection{Relative cells}
Another way to phrase the graded construction of the last paragraph is that we
have defined an operad structures on the relative cells
$CC_*(A,A\setminus\Arc)$.

\subsection{Extending the operadic structure on the cells of $\A^s_{g,r}$}
\label{extasec}
 We can extend the operad structure on
$CC_*(A,A\setminus\Arc)$ to a partial operad structure on
$CC_*(A)$. Recall that $\A_{g,r}^s$ is a CW complex whose cells are
indexed by graphs with possibly isolated vertices. We have to
define the gluings for the boundaries with isolated vertices. For
two arc graphs $(\G,F,[i]),(\G',F',[i'])$ we define $\G\circ_i
\G'$ to be given by the induced operation from the topological
level if $\G$ has arcs incident to the boundary $i$ and $\G'$ has
arcs incident to the boundary $0$.

There are more extensions which are natural. If either $\G$ has no
arcs incident to the boundary $i$ or $\G'$ has no arcs incident to
the boundary $0$ or both, then we define the gluing to be the cell
labelled by the graph obtained from $\G$ and $\G'$ by deleting all
the edges of $\G'$ incident to the $i$th boundary and all  the
edges incident to $\G$ at the vertex corresponding to the $0$th
boundary and the vertices at these boundaries. This graph is
considered to lie on the surface obtained by gluing the boundary
$0$ of $F$ to the boundary $i$ of $F'$.

\subsubsection{An operad structure on $CC_*(\A^s_{g,r})$}
We obtain an operad structure on the Abelian groups
$CC_*(\A^s_{g,r})$ by taking the above definition in the cases of
both boundaries being hit and
 both boundaries not being hit and setting to zero all other products.
  This gluing prescription together with the $\Snn$ action
  by permutations on the labels indeed gives a cyclic operad structure.
We wish to point out that the gluing of empty to empty actually
raises the dimension of the cells by one. This effect is due to
the $\mathbb{R}_{>0}$ action by scaling, see the discussion of
de-projectivized arcs below in \S \ref{extdsec}.

\begin{prop}
The gluings on the spaces $CC_*(\A^s_{g,r})$
defined above together with the action of the groups on the labels
define a cyclic operad structure.
This cyclic operad structure descends to homology inducing
 a structure of a
cyclic operad on $H_*(\A)(n):=H_*(\amalg_{s,g} \A^s_{g,n+1})$.
\end{prop}

\begin{proof}
For the first part only the associativity needs to be checked, which
is straightforward. For the descent to homology,
first notice
that for any two classes $\a\in CC_*(\A)(n)$, $\b\in CC_*(\A)(m)$:
\begin{equation}\label{dgeq}
d(\a\circ_i\b)=d\a\circ_i\b \pm\a\circ_i \b\pm \a\circ_id\b \pm
d\a\circ_id\b \end{equation} where the signs are the natural signs
induced by the orientations of the cells. The only case which is
not straightforward is the case in which the differential deletes
the sole edge incident to one of the  boundaries involved in the
gluing. In this case the composition yields zero, unless there is
also only a sole edge incident to the other boundary involved in
the gluing. Then the gluing $d\a\circ_i d\b$ does not yield zero
in the complex, since the dimension goes up by one when gluing two
boundary components which do not get hit --- as mentioned above.
In this special case, the result is the graph obtained from
$\a\circ_i \b$ in which the sole edge intersecting the separating
curve which is the image of the glued boundaries is deleted. On
the other hand, this is the only summand of $d(\a\circ_i\b)$ that
is not a summand of $d\a\circ_i\b \pm\a\circ_i \b$.  Now for the
descent property, choose any two classes $a\in H_*(\A)(n)$, $b\in
H_*(\A)(m)$. Let $[\a]=a$ and $\b=[b]$ be two representatives
$d\a=d\b=0$. We wish to set $a\circ_ib :=[\a\circ_i\b]$. Due to
the relation (\ref{dgeq}), we see that indeed $[\a'\circ_i
\b']=[\a\circ_i\b]$, for any other representatives $\a',\b':[\a']=[\a],[\b']=[\b]$. With this definition, the cyclic structure
and associativity are clear.

\end{proof}

Examples of the formula (\ref{dgeq}) can be
read off from Figures \ref{boundprobs3} and \ref{boundprobs}.

\begin{rmk}
We wish to point out that in this setting $\GrOCArc$ presents
itself as a mock-$dg$ operad in the sense that the left and right
hand side of (\ref{dgeq}) agree even when restricted to
$\GrOCArc$.
\end{rmk}

\begin{rmk}
We will see below \S\ref{dgextend} that in a suitable restricted case the
mixed gluing can also be used to augment the chain level gluings
to a $dg$ structure.
\end{rmk}

\section{Relations to Moduli Space and other known operads}

\subsection{Definition of Subspaces, Suboperads and $\Darc$.}
\subsubsection{Suboperads}
\label{subspaces} We would like to recall and introduce the
following notation for subspaces.\\

\begin{tabular}{l|l}
Subspace&Condition\\
\hline\\[-2mm]
 $\Arc_g^s(n)\subset \A_{g,n+1}^s$ & $\del:F_{\Gamma}\rightarrow
V_{\Gamma}$ is surjective. This means\\
&  that each boundary gets hit
by an arc.\\

$\Arc^s_{\#g}(n)\subset \Arc_{g}^s(n)$ &the arcs are
quasi-filling. This means that \\& complementary regions are
polygons \\&or once punctured
polygons.\\

$\Tree(n) \subset \Arc(n)$ &  all arcs run only
from boundary $0$ to some\\& boundary $i\neq 0$.\\

 $\Lintree(n)\subset \Tree(n)$ &  the linear order
of  the arcs at the boundary $0$\\& is anti-compatible with the
linear order at each\\& boundary. I.e.\ if $\prec_i$ is the linear
order at $i$ then\\&  $f\prec_i f'$ is equivalent to
$\imath(f')\prec_0
\imath(f)$.\\

$\Corol$&exactly one arc for each boundary $i\neq 0$\\& which runs
to boundary $0$.\\
\end{tabular}\\[2mm]

We will use the subscript $cp$ to signify $g=s=0$: Explicitly
$$\Arc_{cp}(n):=\Arc_0^0(n), \quad
 \Tree_{cp}(n):=\Tree(n) \cap  \Arc_{cp}(n)$$
$$ \Lintree_{cp}(n):=\Lintree(n) \cap \Arc_{cp}(n),
\quad \Corol_{cp}(n)=\Corol(n)\cap \Arc_{cp}(n)$$

\begin{rmk}
One can also define the natural linear order of the flags at the
boundary zero to the one which is opposite to the linear order
induced by the orientation. This convention is in line with the
usual cobordism point of view used in \cite{KLP}. In this case the
condition for $\Lintree$ is the {\em compatibility} of the orders.
\end{rmk}

\subsubsection{De-projectivized arcs}
\label{darcs}
\begin{df}
 Let $\Darc_{g,r}^s := \Arc_{g,r}^s\times \mathbb{R}_{>0}$.
\end{df}
 The elements of $\Darc$ are graphs on surfaces with a
metric, i.e.\ a function  $w:E_{\Gamma}\rightarrow
\mathbb{R}_{>0}$. Furthermore $\Darc$ is a cyclic operad
equivalent to $\Arc$ \cite{KLP}. The operad structure on $\Darc$
is given as follows. Let $\a,\a'$ be elements of $\Darc$, if the
total weight at the boundary $i$ of $\a$ is $\lambda$ and the
total weight at the boundary $0$ of $\a'$ is $\mu$, then first
scale the metric $w$ of $\a$ to $\mu w$ and likewise scale the
metric $w'$ of $\a'$ to $\lambda w'$ and afterwards glue as above.

Any subspace $\mathcal{S}$ of the list above defines a suboperad
$\mathcal{DS}:=\mathcal{S}\times \mathbb{R}_{>0}$ of $\Darc$ which
is equivalent to $\mathcal{S}$.

In the above notation one always has isomorphisms of operads
$\mathcal{DS}/\mathbb{R}_{>0}\simeq \mathcal{S}$ where
$\mathbb{R}_{>0}$ acts by scaling on the right factor
$\mathbb{R}_{>0}$ of $\Darc$.

\begin{lem}
$\Tree,\Lintree$ and $\Corol$ as well as their restrictions to
$g=s=0$ are suboperads (not cyclic) of the cyclic operad $\Arc$.
The same holds for their versions in $\Darc$ defined above.
\end{lem}
\begin{proof}
Straightforward, see also \cite{KLP,cact}.
\end{proof}

\subsection{Relations of subspaces and operads to known operads}
It will be convenient for the reader to list the known
equivalences of sub-operads

\begin{tabular}{l|ll|ll}
Suboperad& \multicolumn{2}{l|}{isomorphic operad}&\multicolumn{2}{l}{equivalent operad}\\
\hline&&\\[-2mm]
$\mathcal{D}\Tree_{cp}$&$\Cacti$ &\cite{KLP}&$fD_2$ & \cite{cact}\\
$\mathcal{D}\Lintree_{cp}$&$\Cact$
\settowidth{\ghost}{$i$}\makebox[1\ghost]{}
 &\cite{KLP}&\settowidth{\ghost}{$f$}\makebox[1\ghost]{}$D_2$ &\cite{cact}\\
 $\mathcal{D}\Corol_{cp}$&$\SCC$&\cite{cact}&\settowidth{\ghost}{$f$}\makebox[1\ghost]{}$A_{\infty}$&\cite{cact}\\
\end{tabular}

The operads in the third column are the familiar ones, that is
 $D_2$ is the $E_2$ operad of little discs, $A_{\infty}$ is the
$E_1$ operad of little intervals and $fD_2$ is the framed little discs operad.
The reader unfamiliar with the second column can take the first column
as an equivalent definition or consult \cite{cact}. Here $\SCC$
is the suboperad of cacti with only one vertex.

Let $M^{1^{n+1}}_{g,n+1}$ be the moduli space of curves of genus $g$
with $n$ punctures and one tangent vector at each puncture.
Below, we will show that:

\begin{tabular}{l|l}
& is isomorphic to\\
\hline & \\[-2mm]
 $\mathcal{D}\Arc_{\#,g}^0(n)$&$M^{1^{n+1}}_{g,n+1}$ \\
$\GrOCArcno$&$\Rib$\\
\end{tabular}

We will first show that the first line is an isomorphism
on the level of spaces, and the second on the level of free Abelian groups.
As we show below, the collection $\mathcal{D}\Arc_{\#}^0(n)$
forms a rational operad. This induces the structure of a rational operad on $\Mngn$.
We also show that $\GrOCArcno$ carries an operad structure,
which is induced by the operad structure of $\Arc$. This operad structure
then carries over to $\Rib$ and hence gives an operad structure on a cell model
of $M^{1^{n+1}}_{g,n+1}$, see \S\ref{cellops} and \S\ref{subspaces}.

\begin{rmk}
The respective quotients by the scaling action or $\Rp$  give rise
to equivalent operads.
\end{rmk}

\begin{rmk}
The inclusion of $\Tree_{cp}\subset \Arcn$ thus gives an  $BV_{\infty}$ (BV up to homotopy) structure to a cell model of moduli which includes an
$A_{\infty}$ structure.

\end{rmk}

\subsubsection{Extended gluing}
\label{extdsec}
 We can extend the operad structure in $\Darc$ on
the spaces $\DA^s_{g,r}:=\A_{g,r}^s\times R_{>0}$ as follows. If
both the boundaries that are to be glued are hit, then glue as
above. If none of the two boundaries to be glued is hit, define
the composition to be given by the image of the two arc graphs on
the glued surface minus the two vertices corresponding to the two
boundaries that are glued together. Lastly if one of the
boundaries is empty, then we delete the vertices corresponding to
the two boundaries that are to be glued along with all the edges
incident to these vertices. The latter operation is not in general
associative, but it can be shown that we obtain an operad up to
homotopy or a quasi--operad (see \cite{cact} of
Definition \ref{quasiopdef} below). These operations
induce the extended chain operations on the space
 $CC_*(A)$ which are discussed in \S\ref{extasec}.
The operations of gluing boundaries which are hit to boundaries
which are not hit still only yield a partial operad structure,
however, since we have to take care that there is at least one arc
left. To obtain an unrestricted operad structure one has to allow
graphs without any edges. We can include the empty arc family as a
point in $\DA^s_{g,n}$ as the image of the origin in ${\mathbb
R}_{\geq 0}^n$, that is the family whose weights are all zero. We
define the space $\dot \DA^s_{g,r}$, to be $\dot
\DA^s_{g,r}:=\A_{g,r}^s\cup \{\oslash\}$ where the topology is
defined by considering $\oslash$ as the limit in which
 all weights of the edges go to zero. This space is obviously
contractible to $\oslash$ by scaling all the weights on the edges
of a given graph homogeneously down to zero. Summing up, we
obtain:

\begin{lem}
There  are contractible spaces $\dot \DA$ which form a quasi--operad whose structure maps induce the operad structure of $\A$ and $\Arc$.
\end{lem}

\qed

\subsection{The duality between quasi-filling arcs and ribbons graphs}
\begin{figure}
\epsfbox{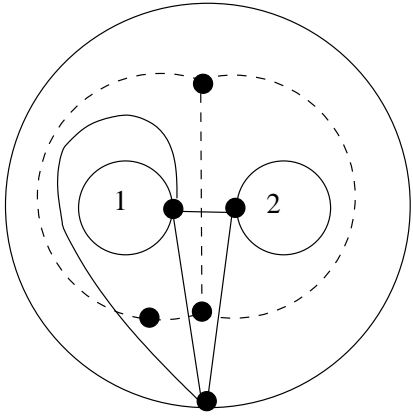}
\caption{\label{dualgraph}
An example of an arc graph (solid) and its dual (dashed).}
\end{figure}

In \cite{KLP} we defined a map called $\mathcal{L}oop$ which is
the suitable notion of a dual graph for a graph on a surface. This
map  uses an interpretation of the graph as a partially measured
foliation. If one restricts to the subspace $\Arcn$ though, this
map has a simpler purely combinatorial description. This
description will be enough for our purposes here, but we would
like to emphasize that this description is only valid on the
subspace $\Arcn$ and cannot be generalized to the whole of $\Arc$
unlike the map $\Loop$.  Figure \ref{dualgraph} contains an
example of an arc graph and its dual.

\subsubsection{The dual graph.} Informally the dual graph of an
element in $\Arcn$ is given as follows. The vertices are the
complementary regions.  Two vertices are joined by an edge if the
complementary regions border the same arc. Due to the orientation
of the surface this graph is actually a ribbon graph via the
induced cyclic order. Moreover the marked points on the boundary
make this graph into a marked ribbon graph. A more precise formal
definition is given in the next few paragraphs.

\subsubsection{Polygons and $\Arcn$}
By definition, in $\Arc_{\#}$ the complementary regions are
$k$-gons or once punctured $k$-gons. Let
$Poly(F,\Gamma,\overline{[i]})$ be the set of these polygons and
let $Sides(F,\Gamma,\overline{[i]})$ be the disjoint union of sets
of sides of the polygons. We define
$\del_{poly}:Sides(F,\Gamma,\overline{[i]})\rightarrow
Poly(F,\Gamma,\overline{[i]})$ to be the map which associates to a
side $s$ of a polygon $p$ the polygon $p$. The sides are either
given by arcs or the boundaries. We define the map
$\lab:Sides(F,\Gamma,\overline{[i]}) \rightarrow E_{\Gamma}
\bigcup V_{\Gamma}$ that associates the appropriate label. Notice
that for $e\in E_{\Gamma};|\lab^{-1}(e)|=2$ and for $v \in
V_{\Gamma}:|\lab^{-1}(v)|=1$.
 Thus there is a fixed point free involution
$\imath_{side}$ on the set $\lab^{-1}(E_{\Gamma})$ of sides of the
polygons marked by arcs which maps one side to the unique second
side carrying the same label. This in turn defines an involution
$\imath$ of pairs $(p,s)$ of a polygon together with a side in
$\lab^{-1} (E_{\Gamma})$ by mapping $s$ to $\imath_{side}(s)$ and
taking the polygon $p$ to the polygon $p':=\del_{poly}(\imath(s))$
of which $\imath_{side}(s)$ is a side. Although $p$ and $p'$ might
coincide the sides will differ making the involution $\imath$
fixed point free.

\subsubsection{The dual graph of an element of $\Arc$} For an
element $\a=(F,\Gamma,\overline{[i]},w) \in \Arcn(n)$ we define
the dual graph
 to be the marked ribbon graph with a projective metric
 $\hat \G(\a):=(\hat \Gamma,\hat w,\mk)$ and $\Zz$ marking on the vertices
$\punc$ which is defined
 as follows. The vertices of $\hat\Gamma$ are the complementary
regions of the arc graph (i.e.\ the polygons) and the map $\punc$
associates to a vertex the number $1$ if the complementary region
is punctured and $0$ if it is not. The flags of the graph  are the
pairs $(p,s)$ of a polygon (vertex) together with a side of this
polygon marked by an arc ($s\in \lab^{-1}(E_{\Gamma})$). The map
$\del$ is defined by $\del((p,s))=p$ and the involution
$\imath((p,s)):=(\del_{poly}(\imath_{side}(s)),\imath_{side}(s))$.

Each polygonal complementary region is oriented by the orientation
induced by the surface, so that the sides of each polygon and thus
the flags of $\hat \Gamma$ at a given vertex $p$ have a natural
induced cyclic order making $\hat \Gamma$ into a ribbon graph.

 Notice that there is a one-one correspondence between edges of
the dual graph and edges of $\Gamma$. This is given by associating
to each edge $\{(p,s),\imath(p,s)\}$ the edge corresponding to the
arc $\lab(s)$.

We define a projective metric $\hat w$ on this graph by associating
to each edge $\{(p,s),\imath(p,s)\}$ the weight  of the arc
labelling the side $s$:  $\hat
w(\{(p,s),\imath(p,s)\}):=w(\lab(s))$, where $w$ is the projective
metric on the arc graph.

To define the marking of the ribbon graph, we first notice that the
cycles of $\hat \Gamma$ correspond to the boundary components of the
surface $F$. Let $c_k$ be the cycle of the boundary component
labelled by $k$. The $k$-th boundary component then lies in a unique
polygon $p=\del_{poly}(\lab^{-1}(k))$. Let $\prec_p$ be the cyclic
order on the set of sides of p, $\del^{-1}(p)$. Let $s_k$ be the
side corresponding to the boundary and let $\Cyc(s_k)$ the element
following $s_k$ in $\prec_p$. We define $\mk(c_k):=(p,\Cyc(s_k))$.

\begin{rmk}

The space $\Arcno$ corresponds graphs with $\punc \equiv0$ and we
will omit this function $\punc$ for these graphs.
\end{rmk}

The above map will suffice for the purposes of this paper. For the
general theory and the reader acquainted with the constructions of
\cite{KLP}, the following  will be helpful.

\begin{rmk}
For elements in $\Arcno$ the dual graph realizes the map $\Loop$
of \cite{KLP}, i.e.\ for  $(F,\Gamma,\overline{[i]},w) \in
\Arc^0_{\#g}(n), \Loop (F,\Gamma,\overline{[i]},w)=(\hat
\Gamma,ord,\hat w,\mk)$.
\end{rmk}

\subsubsection{From marked weighted ribbon graphs to arc
families}

Given  a marked weighted ribbon graph $\G$, we fix $F=\Sigma(\G)$.
The boundary components of $F$ correspond to the cycles of $\G$
and thus the former are labelled if the latter are. If $c_i$ is a
cycle of $\G$ we denote the corresponding boundary by $\del_i F$.
Let $\hat G$ be the dual graph on the surface of $\G$. This graph
can be constructed as follows. Let $\G\subset F$ be embedded as
the spine. For each edge $e=\{f,\imath (f)\}$ with $f\in c_i$ and
$\imath (f)\in c_j$. Fix the mid-point of each edge $mid(e)$. Now
let $\hat e$ be the arc $mid(e)\times I$ on $F$. This edge is
broken into two flags $\hat f,\mathcal (\hat f)$ by the midpoint
of the interval. Here the flag which is named $\hat f$ is the flag
of $\hat e$ which runs from the midpoint $mid(e)$ to the boundary
$\del_i F$. This defines a map $\hat{}: F_{\G}\rightarrow F_{\hat
\G}$. It then follows that $\widehat{\imath (f)}=\imath(\hat{f})$.
From $\hat \G$ it is easy to obtain an element of $\Arc$ in the
picture where the arcs do not run to the marked point on the
boundary.
 Just mark a point on each
boundary $\del_i F$ such that the linear order of incident edges to
$\del _i F$ is that of the cycle $c_i$. For instance any point which
is slightly before the point of intersection of
$\widehat{\mk(c_i)}$. To obtain a graph on $F$ in the sense of \S
\ref{arcsection}, let $p_i\in \del_i F$ be the endpoint of
$\widehat{\mk(c_i)}$ and choose a homotopy of the flags $\hat{f}$
which fixes the point $mid(e)$ and slides all the endpoints on the
various $\del_i F$ to fixed point $p_i\in \del_i F$ in the direction
opposite to the natural orientation of $\del_i F$. Obviously the
homotopy can be chosen, such that all the arcs will be embedded. The
edges which are now called arcs incident to $p_i$ will then have a
linear order starting with the edge to which $\widehat{\mk(c_i)}$
belongs.

\begin{prop}
\label{graphident} The dual map identifies the space $\DArcno$
with the space of metric marked ribbon graphs $\MRib$. Moreover
this correspondence is on the level of $\mathbb{R}_{>0}$ spaces
where the action on $\MRib$ is by rescaling the metric, i.e.\
$\Arcn$ and $\PRib$ are identified by the dual graph construction.
\end{prop}

\begin{proof}
As seen above we have maps $\phi:\DArcno\rightarrow \MRib$ and
$\psi:\MRib\rightarrow\DArcno$. It is straightforward to check
that $\phi$ and $\psi $ are inverses, since the surface of an
element in $\DArcno$ can uniquely be reconstructed from the
underlying graph.
\end{proof}

\begin{rmk}
If we regard pairs $(\G,\punc)$ then we obtain a map to $\Darcn$
where in the last step we put a puncture at each vertex $v$ which
is has $\punc(v)=1$. There is no operad structure however, since
we cannot guarantee that the number of punctures does not
increase.

 Of course, we could keep track of this and regard the subspace $\Arc_{\#\#}$ of elements
of $\Arc$ whose complementary regions are polygons, or multiply
punctured polygons. Then the function $\punc$ would take values in
$\mathbb{N}$ and we would  obtain analogous statements to the
Proposition \ref{graphident} and Corollary \ref{ribop}.
\end{rmk}

\subsection{Ribbon graphs and Moduli-Space}

The spaces $\Arc_g^s(n)$ are related to moduli spaces of Riemann
surfaces with extra structures. This for instance follows using
Strebel differentials \cite{Strebel} to
relating the $\Arc$ to the moduli space of curves with marked
points and tangent vectors at these points.

\begin{thm}
The space $\Darc_g^0(n)$ and the moduli space
$M_{g,n+1}^{1^{n+1}}$ of $n+1$ punctured Riemann surfaces of genus
$g$ with one tangent vector fixed at puncture each are isomorphic
in the coarse moduli or orbifold sense. Moreover they are isomorphic as
$\mathbb{R}_{>0}$ spaces where the action on $M_{g,n+1}^{1^{n+1}}$
is given by simultaneous rescalings of all tangent vectors, i.e.\
$\Arcno\simeq
M_{g,n+1}^{1^{n+1}}/\mathbb{R}_{>0}:=\mathbb{P}M_{g,n+1}^{1^{n+1}}$.
\end{thm}

\begin{proof}
As shown each $\alpha \in \Darc_g^0(n)$ uniquely corresponds to a
marked ribbon graph with a metric, which is embedded as the spine of
the given surface. Now as usual by gluing in punctured discs to the
boundaries using Strebel differentials, we obtain a punctured
surface. Now, furthermore, we retain the length of the boundary
cycle and a tangent direction. This direction is given by the
direction vertical trajectory which hits the marked point on the
boundary. This data gives rise to a surface with a marked point and
a tangent vector at each boundary, by decomposing $v\in T_{p_i}F$ as
$|v|e_{v}$ with $e_v\in T^1_{p_i}F$. It is clear that this
assignment is a bijection in view of the known characterization of
combinatorial moduli space \cite{Strebel,Koint,
formal,HubMas,Harstab,P3,P2, Vogt}. Moreover this correspondence is
compatible with the orbifold topology or the coarse moduli
structure. The statement about the $\mathbb{R}_{>0}$ action is
obvious from the description.
\end{proof}

\begin{cor}
There is a decomposition of $M_{g,n+1}^{1^{n+1}}=\amalg_{\G\in
\Rib} \dot C(\G)$ with $\dot C(\G)$ open cells indexed by marked
ribbon graphs.
\end{cor}

\begin{rmk}
Alternatively, one could use Penner's formalism \cite{P2} in the
hyperbolic setup to get a proper homotopy equivalence between
$\Arc_{\#}(F)$ and $M(F)/\mathbb{R}_{>0}$  for any bordered
surface $F \neq F^0_{0,2}$. Here $M(F) = [Hyp(F) \amalg (\prod^r_1
\del_i)]/ \sim $ where $\sim$ is the equivalence relation
generated by the push-forward of the metric under
orientation--preserving diffeomorphism and $Hyp(F)$ is the space
of all hyperbolic metrics with geodesic boundary on the surface
$F$ whose boundary components are the $\del_i$.
\end{rmk}

\subsubsection{The Graph Complex and Homology of the Mapping Class Group}
\label{graphcompsection} By using the ribbon graph interpretation
for $\OCArcn(F)$ and its differential we obtain  the graph complex
of marked ribbon graphs.  Furthermore its homology computes the
cohomology of the mapping class group of $F$ by the usual
arguments \cite{formal,P2,P3,Vogt}.

\begin{nota}
We let $\OCArcno$ be the subgroup of $\OCArc$ generated by the
cells corresponding to quasi-filling arc families with no
punctures and write $\GrOCArcno$ for the image of this subgroup on
$\GrOCArc$.
\end{nota}

\begin{df}
The graph complex of marked ribbon graphs is the Hopf algebra
whose primitive elements are connected marked ribbon graphs and
whose product is the disjoint union. Its differential is given by
the sum of contracting edges $d\Gamma=\sum_{e\in E'(\G)} \pm
\Gamma/e$, where $E'(\G)$ is the subset of edges $e$ s.t.\ the
topological type of $\G$ coincides with that of $\G/e$  and the
sign is the usual sign.
\end{df}

\begin{prop} The following two complexes are isomorphic
\begin{itemize}
\item[i)]$(\GrOCArcno(F),d)$ \item[ii)] Graph complex of marked
ribbon graphs
\end{itemize}
and they both compute $H^*(\PMC(F))$, the cohomology of the pure
mapping class group and the spaces $H^*(M_{g,n+1}^{n+1})$.
\end{prop}

\begin{proof}
The differential of the graph complex is the signed sum over those
contractions of edges, which leave the topological type of the
graph (genus and number of boundaries) intact. This is exactly
dual the differential of $\OCArcno$ which deletes the respective
arcs. The fact that they both compute $H^*(\PMC(F))$ is essentially
the Kontsevich-Penner theorem adapted to the case with boundary.
  The proof is a standard application of the techniques of
\cite{formal,P2,P3,Vogt}.

\end{proof}

\section{The cell operad of moduli space}

Notice that the gluing operation of $\Arc$ applied to two elements
of $\Arcno$ need not land in $\Arcno$. Generically, however, that is
inside the top--dimensional cells, two elements in $\Arcno$ do glue
to form an element in $\Arcno$. In order to deal with this
situation, we introduce the following notion.

\begin{df}
A rational topological operad is a collection of topological $\Sn$
modules together with operadic structure maps which only need to
be defined and continuous on a dense subset. These structure maps
are required to satisfy the operad axioms where they are defined.

We will call a rational operad almost topological if the structure
maps of the rational operad can be extended to the whole space in
a possibly non--continuous fashion, such that forgetting the
topology, the induced structure is an operad of sets.
\end{df}

\begin{thm}
\label{fillingopprop} $\Arcno$ is a rational topological operad
and $\GrOCArcno\subset\GrOCArc$ is a cyclic suboperad which is a
$dg$-operad.
\end{thm}

\begin{proof}
We will show that the top--dimensional cells in the composition of
two cells from $\OCArcno$ are also cells of $\OCArcno$. This can be
done with the help of the Euler characteristic. First notice that
the top--dimensionality of the cells implies, that there are no
closed loops and that we are in the generic situation for the gluing
as described in \S\ref{gluingpar}. Notice that by Lemma \ref{Euler}
we only have to show that the Euler characteristic of the arc family
is the Euler characteristic of the glued surface.  Now gluing the
arc families supported on the surfaces $F_1$ and $F_2$, say, we
obtain a graph $\G$ on the glued surfaces which we call $F$. Let $S$
be the separating curve which constitutes the glued boundaries of
the two glued surfaces $F_1$ and $F_2$.
 Now the
Euler characteristic satisfies $\chi(F)=\chi(F_1)+\chi(F_2)$.  Let
$\G_1$ on $F_1$ and $\G_2$ on $F_2$ be the graphs before gluing.
Fix a partition or $(k,l)$-shuffle (cf.\ \S\ref{gluingcomb}) indexing a top-dimensional cell
and let $\G'_1$ and $\G_2'$ be the graphs obtained from $\G_2$ and
$\G_1$ by inserting parallel $k+l$ edges and breaking up the
vertices $0$ and $i$ according to this partition in the gluing
procedure of \S\ref{gluingpar}. See Figure \ref{gluingex1} for an example.
Let $\G'$ be the graph consisting
of $\G'_1$, $\G'_2$ and the curve $S$. It is clear that $\G\subset
\G'$.

The graph $\G$ is obtained from
$\G'$ by erasing the edges belonging to $S$. We claim that
\begin{equation}
\label{*}
 \chi(\G')=\chi(\G_1)+\chi(\G_2)
\end{equation}
and that
\begin{equation} \label{**}
\chi(\G')=\chi(\G).\end{equation}
Assuming these equalities we obtain
$\chi(F)=\chi(F_1)+\chi(F_2)=\chi(\G_1)+\chi(\G_2)=\chi(\G')=\chi(\G)$,
whereby $\chi(F)=\chi(\G)$ and hence by Lemma \ref{Euler} $\G$
is quasi-filling. To validate (\ref{*})  we calculate:
$V(\G')=V(\G_1)+V(\G_2)+(k+l-1)$, $E(\G')=E(\G_1)+E(\G_2)+2(k+l-1)$
and $\C(\G)=\C(\G_1)+\C(\G_2)+k+l-1$ so that
$V(\G')-E(\G')+\C(\G')=V(\G_1)-E(\G_1)+\C(\G_1)+V(\G_2)-E(\G_2)+\C(\G_2)$. To
check (\ref{**}), we first need to make some observations. First,
the curve $S$ breaks up into the  pieces corresponding to the
subintervals of $\mathcal{P}_{\circ}$ of \S \ref{gluingcomb} in
their order
---which we will call intervals--- and  one piece which is between
the first and the last arc ---which we will call the outside arc
of $S$. Secondly, the complementary regions of $\G'$ correspond to
those of $\G_1$, $\G_2$, triangles corresponding to non-recursive
subintervals  and quadrangles corresponding to recursive
subintervals. Notice that at least on one side of an interval is a
piece which has parallel flags, viz.\ a quadrangle or a triangle,
since we are in the maximal dimension and we have a full
partition. Moreover if a subinterval is the boundary between a
polygon and a quadrangle, then on the opposite side of this
quadrangle will be another quadrangle or triangle, since the flags
will remain parallel. Now because of these observation, we see
that removing the intervals does not change the Euler
characteristic. Moreover, we see that no two of the original
polygons of $\G_1$ and $\G_2$
 have been joined when removing the intervals.
 So finally by removing the outside arc of $C$ one
effectively glues two distinct polygons together and hence again
does not change the Euler characteristic.

For the $dg$ part, we notice that in view of the considerations
above if an arc of $\G$ is not one of the glued ones then the
condition  that it is removable ---viz.\ the condition that
the remaining graph is quasi-filling---
is the same  before and after the gluing.
Moreover the same holds true for the arcs which are glued. If a
weight of one of these goes to zero, then this corresponds to a
face of the simplex defined by $\G$, cf.\ \S \ref{gluingcomb}.
If this face is allowed
in $\G$, that is $\G$ after deleting this edge is still quasi--filling,
 the respective limit is allowed by the above in $\G_1$ and
$\G_2$. The same holds true for all the iterations of passing to faces.
\end{proof}

Notice that in the quasi-filling case the limit which appeared in
\S\ref{nodg} as a counter-example to the compatibility of the
operad structure with the differential in the general case is not
valid anymore since the relevant limit is neither allowed for
$\G_1$ nor for $\G_2$ inside $\Arcno$.

\begin{rmk}
The equation (\ref{*}) above also holds in the case that both the
families are twisted at the boundaries which are glued.
\end{rmk}

\begin{cor}
\label{ribop}
$\mathcal{D}\Arc_{\#}$ and $\Mngn$ are a rational operads.
Moreover there is an operad structure on the free Abelian
group generated by ribbon graphs $\Rib$ which is induced via the
identification $\Rib\simeq \OCArcno\simeq \GrOCArcno$.
That is there is a $dg$--operad structure on the graph complex computing
the cohomology of $\Mngn$.
\qed
\end{cor}

\section{Di-operads and PROPs based on the $\Arc$ operad}

\begin{nota}
In all the following the subscript $\#$ will mean that the
condition that arc families under consideration are quasi-filling.

{\sc Convention.} In this paragraph, we will also restrict to the case $s=0$ when we restrict to the
quasi--filling case.
\end{nota}

\subsection{The Di-operad $\Diarc$} We will consider additional
markings for elements of the $\Arc$ operad. The first is a
partition of $\bar n=\In \amalg \Out$. Such a partition is
equivalent to a map $\io:\bar n \rightarrow \Zz$ where $\In
=\io^{-1}(1)$ and $\Out=\io^{-1}(0)$. Let $Part(n)$ be the space
of maps from $\{0,\dots,n\}$ to $\{0,1\}$ and let
$Part(k,l)\subset Part(k+l-1)$ be the subset of functions with
$|\In|=k,|\Out|=l$.

Let $\Z/2\Z[1]=\{\Z/2\Z[1](n)\}$ be
 cyclic operad of spaces built on $\Z/2\Z$.
That is the operad of spaces as defined in
\cite{cact} shifted by $1$.
$\Z/2\Z[1](n)=\Z/2\Z(n+1)=(\Z/2\Z)^{\times n+1}$ where the
indexing set for the Cartesian product is taken to be $\bar n=\{0,
\dots,n\}$ and the action of $\Snn$ is by permutations.

Set $\Diarc(n,m):=\Arc(n+m-1) \times Part(n,m)$. Identifying
$Part(n,m)$ with $\Z/2\Z^{n+m-1}$ and restricting the gluing in
the cyclic $\Arc$ operad to only gluing ``ins'' to ``outs'' and
retaining the in/out designation on the non-glued boundaries one
immediately obtains:

\begin{prop}
Let $\Diarc$ be the collection of $\;\Sn\times \Sm$ modules
$\Diarc(n,m)$ where the action of the symmetric groups is the
action of permuting the labels of the $\In$ and $\Out$ boundaries,
then is  $\Diarc$ a di-operad. It is isomorphic to a partial suboperad of
the direct product of operads $\Arc\times \Z/2\Z[1]$. Furthermore
$\Diarcno$ is a sub-di-operad on the cellular level, i.e.\ the
respective subgroups in $\GrOCArcno(n)\times(\Z/2\Z)^{n+1}$ form a
di-operad.
\end{prop}
\qed

\begin{df}

We define the following sub-spaces (recall that the subscript \#
also implies $s=0$ in this section).

\begin{tabular}{l|l}
Sub-space&Condition\\
\hline
 $\iooarc(n,m) \subset \Diarc(n,m)$&
arcs from input to output boundaries \\&or from output to
output boundaries\\& only. Viz.\ no arcs from input to input.\\
$\ioarc(n,m) \subset \iooarc(n,m)$ & arcs only from input to
output boundaries.\\

$\Sularc(n,m)\subset \iooarc_{\#}(n,m)$& after removing all arcs
which run from out\\& to out the arc family is still quasi-filling.\\

 $\StSularc(n,m)\subset \StSularc(n,m)$& arc families
 such that no two neighboring\\& arcs on an ``out'' boundary
both run to\\& ``in'' boundaries.
\end{tabular}
\end{df}

\subsection{Relation to graphs and Sullivan chord diagrams}
\begin{prop}
The dual graphs of $\ioarc_{\#}(n,m)$ lie in $\ppdigraph(n,m)$,
the dual graphs of $\StSularc(n,m)$ lie in $\StSulchord$, and the
dual graphs of $\Sularc(n,m)$ lie in $\Sulchord(n,m)$. Moreover
the mentioned subspaces of $\ioarc$ are naturally identified with
the relevant subspaces of $\;\PP\Digraph$ of graphs with a
projective metric. Therefore the cells of $\OCDiarc\simeq
\GrOCDiarc$ belonging to these subspaces are exactly indexed by
the graphs of the indicated type.
\end{prop}
\begin{proof}
The claim about the graphs becomes clear by unravelling the dual
graph construction. A dual edge to an arc which runs from in to
out is an edge which is part of an $\In$ and an $\Out$ cycle.
Likewise an arc from out to out yields an edge which belongs to
two  (not necessarily distinct) $\Out$ cycles. The condition that
the out-out edges form trees is equivalent to the fact that
contracting them does not change the genus and the number of
cycles. In other words the contracted graph and the graph define
the same topological surface. Contraction is dual to deletion
hence the condition to be quasi-filling is dual to the condition
of contractibility. The condition on Sullivan Chord diagrams in
the strict sense is that the $\In$ cycles are disjointly embedded.
This means that no two $\In$ cycles share a vertex, i.e.\ there is
at least one out-out edge between them. The dual to this condition
is the one that is stated. Now the rest of the statements directly
follow.
\end{proof}

\begin{rmk}
Just like $\Arcn$, the subspaces $\StSularc$
 and $\Sularc(n,m)$ are
not stable under composition, since the condition
 of non-topology changing contractibility of the out-out arcs
 is not stable under the composition.
\end{rmk}

\begin{rmk}
It is clear that $\Sularc(n,m)$ retracts to its subspace
$\ioarc_{\#}$ by the homotopy that homogeneously scales all weights
on the arcs from ``outs'' to ``outs'' to zero.

Just as the respective graphs, see Remark \ref{weaksul},
$\ioarc_{\#}(n,m)$ and $\StSularc(n,m)$ are weakly homotopy
equivalent.
\end{rmk}

\begin{prop} The collections of $\;\Sn\times \Sm$ modules,
 $\iooarc(n,m)$, $\ioarc$ form di-operads.

The subspaces $\iooarc_{\#}$ and $\ioarc_{\#}$ are rational
sub--dioperads which induce $dg$--di--operads on
 the graded open-cell operad level.
Hence there is an induced $dg$ di-operad structure on the respective
graphs.
\end{prop}

\begin{proof}
On the topological level, we only have
to show that the gluing preserves the subspaces.
In the operadic composition, the arcs/bands are matched or split
and then matched. In both cases an arc which runs from in to out
will be continued with an arc from in to out and thus will run
from in to out. In the case of $\iooarc$ an arc which runs from
out to out might be matched with an arc running from in to out.
The resulting arc will again run from out to out. Therefore the
di-operad structure follows from the operad structure of $\Arc$.

Lastly, the claims about the open cell level all follow from
Proposition \ref{fillingopprop}, the remarks above and the fact
that $\GrOCArcno$ is a dg-operad and hence the respective boundary
limits before and after gluing coincide.
\end{proof}

\subsubsection{Standard In/Out markings}
\label{standardmarkings} There are no natural partitions into in
and output boundaries for general elements of $\Arc$, except all
inputs or all outputs; it is after all a cyclic operad. If one has
a family, however, which has a partition of the boundary
$S_1\amalg S_2=\{0,\dots, n\}$ such that the arcs of this element
only run between $S_1$ and $S_2$ then one has a $\Zz$ choice of
calling $S_1$ either $\In$ or $\Out$. In this case, we will fix that the set
containing $0$ will be called $\Out$. This establishes an
identification of $\ioarc$ with a partial suboperad of $\Arc$. Also there
is a partition of the boundary $S_1\amalg S_2=\{0,\dots, n\}$ such
that the arcs of this element only run between $S_1$ and $S_2$ or
between $S_2$ and $S_2$ and there is at least one such arc, we
set $S_2=\Out$ and hence identify $\iooarc$ with a partial suboperad of
$\Arc$.

The subspace $\Tree$ is also of the form described above. Hence
the boundaries have a standard $\Zz$ marking with $\io(0)=0$ and
$\io(j)=1$ for $j\neq 0$. This identifies $\Tree\subset \ioarc$ as
a sub-di-operad.

\subsection{Operads from Arc families with angle markings}

There are basically three approaches to generalize the actions of
\cite{del,cyclic}. The main observation is that in those actions
not all boundaries were treated equally. The boundary $0$ played
the special role of an output. The first approach to this
non--symmetric situation
 is to stay in the PROP setting
by  explicitly marking the boundaries as $\In$ or $\Out$. This of
course breaks the cyclic operad structure. The second approach,
which we explain below is to put an additional angle marking on
the operad which preserves the cyclic structure. The third way to
proceed is to partially merge these two approaches, by keeping the
$\In/\Out$ distinction, and using this to define an angle marking.
The angle marking is intimately linked to the operations defined
by these graphs \cite{hoch2}.

\begin{df}
An angle marking for an element
$\alpha=(F,\Gamma,\overline{[i]},w)\in \A_{g,r}^s$ is an angle
marking $\amark$ for $\G$. This is clearly $PMC$ invariant data.
\end{df}
We define
\begin{equation}
\A_{g,r}^{\angle \; s}:=\{(\alpha,\amark)|\alpha
=(F,\Gamma,\overline{[i]},w)\in
\A_{g,r}^s,\amark:\angle_{\G}\rightarrow \Zz\}
\end{equation}
to be the CW complex obtained in the analogous fashion to
$\A_{g,r}^s$ where now  the differential is given by deleting edges
from $(\G(\a),\amark)$ as an angle marked ribbon graph (cf.\
\S\ref{deledge}). We will also use notation analogous to the
notation \ref{ribnota}.

\begin{nota}
For all the subspaces $\mathcal{S}$ listed in \S \ref{subspaces},
we denote the corresponding subspaces as $\mathcal{S}^{\angle}$.
E.g.\ $\Anarc$.
\end{nota}

\begin{rmk}
$\Anarc(n)$ is filtered again filtered by the number of edges
minus one, that is dimension of the cell corresponding to the graph.
Now $|\angle(\a)|=|F(\a)|=2|E(\a)|$.
Let ${\Anarc}^{\
k}(n)$ be the subspace of graphs with  $k+1$ edges then
there is an induced exhaustive filtration: $$\dots \subset {\Anarc}^{\leq
k}(n) \subset{\Anarc}^{\leq k+1}(n) \subset\dots$$
 We will use the
identification ${\Anarc}^k(n):={\Anarc}^{\leq k}(n)/{\Anarc}^{\leq
k-1}(n)$.
\end{rmk}

\begin{rmk}
$\Anarc(n)$ has a graded open cell decomposition
\begin{equation}
\Anarc(n)=\amalg_{(\a,\amark):  \a\in \Arc(n),
\amark:\angle_{\G}(\a)\rightarrow \Zz} \dot C(\G(\a))
\end{equation}
Furthermore using the angle--edge--correspondence
\begin{equation}
{\Anarc_g}^{k}(n)=\amalg_{(\a):  \a\in \Arc_g^k(n)} \dot
C(\G)\times (\Z/2\Z)^{2(k+1)}
\end{equation}
\end{rmk}

\begin{nota}
We let $\Anarcn(n)$ be the subspaces of elements whose arc
families are quasi-filling.
\end{nota}

\begin{rmk}
$\Anarcn(n)$ has a graded open cell decomposition
\begin{equation}
\Anarcn(n)=\amalg_{(\G,\amark) \G\in \MRib(n)} \dot C(\G(\a))
\end{equation}
Furthermore
\begin{equation}
{\Anarcn}_g^{k}(n)=\amalg_{\G\in \MRib_g^k(n)} \dot C(\G)\times
(\Z/2\Z)^{2(k+1)}
\end{equation}
\end{rmk}

\subsubsection{The operadic compositions for $\Anarc$}
\label{anarcop}

 Given $\alpha =(F,\Gamma,\overline{[i]},w)\in \Arc(n)$ and $\beta
=(F',\Gamma',\overline{[i']},w')\in \Arc(m)$ let
$\alpha\circ_i\beta=(F'',\G'',[i''],w''])$. We remark that using
the gluing formalism of \S\ref{gluingpar} the angles of $\G''$
were either formerly angles of $\G$ and $\G'$ which we called
non-split or split angles. On one hand each split angle
corresponds to a triangle in the gluing process before removing an
interval in the notation of Proposition \ref{fillingopprop}. The
interval on the other hand corresponds to an angle of $\G'$ or
$\G''$ which is given by the two flags of the interval on the
opposite triangle, we call this angle the opposite angle of the
split angle. If there is a quadrangle on the opposite side, we
continue the process until we hit a polygon which is not a quadrangle
of $\G$ of $\G'$.  In
both cases, when removing the intervals the identification of the
corresponding flags associates a unique opposite angle to each
split angle.

\begin{figure}
\epsfxsize =\textwidth
\epsfbox{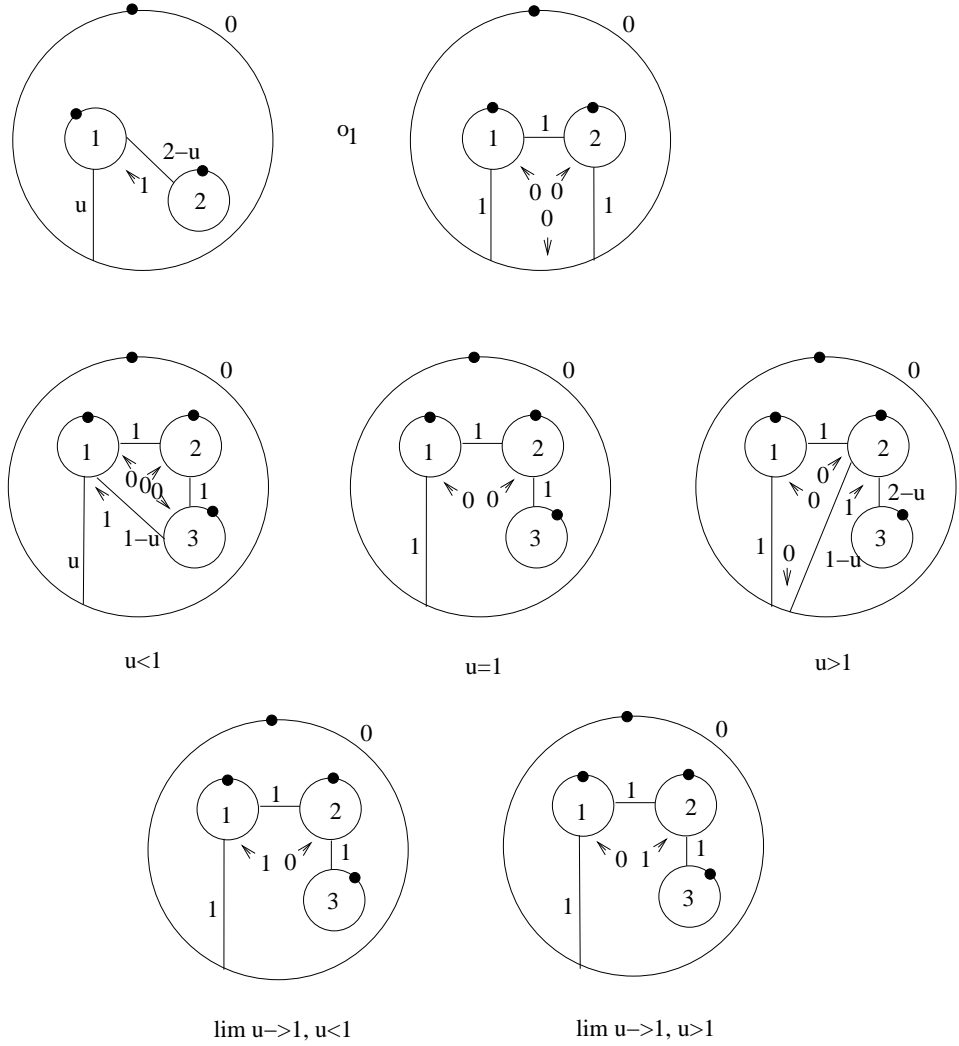}
\caption{\label{limfig}{
An example showing that $\Anarc$ is only a rational topological operad}
}
\end{figure}

Given the angle markings $\amark:\G\rightarrow \Zz$ and
${\amark}':\G'\rightarrow \Zz$ we define
$\amark\circ_i{\amark}':\G''\rightarrow \Zz$ as follows:
\begin{itemize}
\item[i)] If $\theta\in \angle_{\G''}$ is not split then the label
is retained.

\item[ii)] If $\theta\in \angle_{\G''}$ is  split then $\theta$
will be labelled by the label of the angle opposite the split angle.
(See Figure \ref{limfig} examples.)
\end{itemize}

\begin{df}
In the above notation we define
$$
(\alpha,\amark)\circ (\beta,{\amark}'):=(\alpha\circ_i \beta,
\amark\circ_i {\amark}').
$$
\end{df}

\begin{rmk}
Using these gluings, we do not get a topological operad structure
on $\Anarc$, due to the fact that the boundaries are not behaved
well with the natural differential and hence the gluings are not
continuous. (see the example in Figure \ref{limfig}). However, the
gluings are defined everywhere and are continuous up to a
codimension one set.
\end{rmk}

\begin{prop}
The above operations $\circ_i$ together with the $\Sn$ actions
acting by permuting the labels imbue $\Anarc$ with the structure
of a cyclic rational operad that is almost topological. This
operad structure respects the filtration $\leq k$ and hence
induces an operad structure on the associated graded of the open
cell decomposition. Lastly, this operad structure induces a cyclic
operad structure on graded open cell level for  $\Anarcn$ and
hence on the set of
 Abelian groups $\Anrib$ which are isomorphic to the graded open cell
 decomposition of $\Anarcn$.
\end{prop}

\begin{proof}
The fact that the $\circ_i$ yield an operad of sets follows from
the associativity of the operad structure on $\Arc$ and the
associativity of the marking function under the composition. It is
clear that the operad structure is continuous on the interior of
the cells and possibly discontinuous only on the boundaries which
are at least codimension one. Thus we obtain a rational operad
structure that is almost topological. Since the operad structure
on $\Arc$ respects the filtration, so does the operad structure on
$\Anarc$. The other facts are now straightforward.
\end{proof}

\subsubsection{De-projectivized version and extended gluings}
As was the case for $\A$ and $\Arc$ there are straightforward
generalizations of the gluing operations to $\Aang\times \Rp$
and to the analogue of $\dot{\DA}$ and moreover extensions on these
sets involving gluing to empty boundaries. The constructions are
{\em mutatis mutandis} the same as in \S\ref{extasec} and
\S\ref{extdsec}.

\subsubsection{Standard markings and suboperads}
\label{markingpar}

Notice that for an element  $\a\in\Arc,
\alpha=(F,\Gamma,\overline{[i]},w)$ there are two types of angles.
The angles $(f,\Cyc(f))$ in which $f <  (\Cyc(f))$ in the natural
{\em linear order} $<$ at the boundary of $f$
---we will call these angles
the inner angles--- and the angles $(f,\Cyc(f))$ in which $f>
(\Cyc(f))$ which we will call outer angles. There is exactly one
outer angle at each boundary.

There are several embeddings  $\Arc\subset \Anarc$. Each one is
given by choosing a function $\amark$. Three of these choices are
rather canonical. The first is by a constant marking of all angles
by $0$, the second is given by the constant angle marking $1$.
Lastly, one can mark all inner angles by $0$ and the outer angles
by $1$.

All these  markings are interesting and embed $\Arc$ as a cyclic
operad; that
is the image is truly a cyclic operad and not just a rational operad.
The second marking will play a special role for us in \cite{hoch2}
as it leads to a connection with the $\Xi_2$ operad of \cite{MScosimp}.
So for $\a\in \A$ we define the standard marking to be given by:
\begin{equation}
\label{arcmarki} \amark(\theta)\equiv 1
\end{equation}

\begin{cav}
For the non-cyclic suboperads $\Tree$ and $\Lintree$ there are two
standard markings one as suboperads of $\Arc$, and the second as
sub-di-operads of $\Diarc$ as defined in the next paragraph.
\end{cav}

The standard marking $\amark$ for $\Diarc$ is
given by

\begin{equation}
\label{diarcmark}
 \amark (\theta)=\begin{cases}1& \text {if $\theta$ is an
outer
angle}\\
1&\text {if $\theta$ is an inner angle belonging to an $\In$ boundary }\\
0&\text {if $\theta$ is an inner angle belonging to an $\Out$
boundary }
\end{cases}
\end{equation}

This also gives standard markings for $\Tree$ and $\Lintree$, when
considering the boundary $0$ as the only $\Out$ boundary and the
other boundaries as $\In$ boundaries.

\begin{nota}
Sometimes it will be necessary to distinguish between an arc graph
and the arc graph with standard angle markings. To facilitate
this distinction, for an arc graph $\a$ we write $\a^{\angle}$ for
the arc graph with the standard marking. Here the standard marking
for $\Arc$ is defined by (\ref{arcmarki})  and
the standard marking for $\Diarc$ is defined by (\ref{diarcmark}).
\end{nota}

\begin{lem}
The spaces $\Arc$, its suboperads $\Tree$, $\Lintree$ and their
restrictions $cp$ with the above marking (\ref{arcmarki}) or any of the other markings discussed in \S\ref{markingpar}
 are suboperads of $\Anarc$.

The subspace  and $\ioarc$ (and hence $\Tree$ and $\Lintree$
viewed as in \S \ref{standardmarkings} with the convention
(\ref{diarcmark})), are partial suboperads.  Here the partial
compositions are dictated by the marking $\amark$. A boundary is
an $\In$ boundary if one and hence all of its inner angles are marked
by $1$ and a boundary is an out boundary if one and hence all of
its inner angles are marked by $0$. The partial operad structure
is given by restricting the operadic composition to connecting $\In$
to $\Out$ boundaries only.

\end{lem}

\begin{proof}
The claim for $\Arc$ and hence its suboperads, is clear, since all
markings will always be constantly $1$ in the standard marked
case. In the other cases either we have a constant marking or only
the outer angles are marked by $1$. These latter marking is
invariant by the definition of deleting edges in a tree with angle
markings as defined in \S\ref{deledge}.

For $\ioarc$ first notice that the outer angles are always
non-split, so they retain their marking by $1$. Thus in the case
of $\ioarc$, we only have to look at split angles. In the case an
angle is split then the opposite angle is  on a boundary with the
same $\In/\Out$ marking due to the definition of $\ioarc$.
\end{proof}

\begin{prop} The suboperad
 $\Arc(n) \subset \Anarc(n)$ considered embedded via (\ref{arcmarki}) or
 any of the other embeddings of \ref{standardmarkings} is a topological sub--operad. And
 $\ioarc(n,m) \subset \Anarc(n+m-1)$ is a cyclic partial
 topological operad that defines a di-operad in any of the above embeddings.
\end{prop}
\begin{proof}
The only problems that can arise are on the boundary. The
compatibility of the different limit for these sub--spaces follows
from the definition of deleting edges in angle marked graphs. In
particular, for $\Arc$ this observation is trivial, since the angle
marking is constant and stays constant. For $\ioarc$ notice that
since the taking a limit into the boundary by deleting an arc
corresponds to deleting an edge in the dual graph.  The rule for the
boundary marking is given by $\bar a \bar b$. If there is an outer
angle involved, then we have no problem since $\overline {\bar 1\bar
a}=\overline 0 =1$ and the merged angle is again an outer angle. On
inner angles, we always have $a\equiv b$, so that $\overline{\bar a
\bar b} = a$. In other words, the limit obtained by deleting edges
from the angle marked graphs reproduces the standard marking. The
partial operad structure is cyclic and allows to glue one ``in'' to
one ``out'', so by definition it gives rise to a di-operad
structure.
\end{proof}

\begin{cor}
The PL-chain level operads of $\Arc$, $\ioarc$ thought of as
topological suboperads of the chains of $\Anarc$ are dg-chain
operads.
\end{cor}

Just as in the case of $\Arc$ there are analogues of the cell complexes
and quasi--filling sub--complexes:

\begin{prop}
The relative chains $\OCAnarc$
form an operad. The quasi-filling analogues $\Arcn$ and
$\ioarc_{\#}$ and also the associated graded relative chain
complexes of these spaces form $dg$ operads.
\end{prop}

\begin{proof}
 These
statements follows directly from the operad structure of $\Anarc$.
Namely, we can write the degree $k$ component of $\Anarc(n)$ as
${\Anarc}^k(n)=\Arc(n)^k\times (\Z/2\Z)^{2(k+1)}$. We also know that the
composition respects the filtration of $\Anarc$ and the factor
$(\Z/2\Z)^{2(k+1)}$ can be thought of as indexing the cells by discrete
data. So  the statements follow from the fact that the analogous
statements are true for $\Arc$.
\end{proof}

\subsection{Local scaling and ``up to homotopy'' structures}

The di-operad $\Diarc$ falls short of being a PROP, since we can
only guarantee the compatibility of the weights for one boundary by
using suitable representatives in the projective class.  If we are
willing to
relax the associativity up to homotopy, we can,
however, achieve a quasi-PROP structure on $\iooarc$.
 For this we will introduce a new composition
which is given by a local scaling, since globally, we can in general
only scale to match one boundary. Now we will scale only the arcs
incident to the two boundaries which are to be glued. This will
destroy the associativity, but the associativity does hold up to
homotopy so that there will be an honest associative structure on
the homology level and in some situations, with a careful choice of
chains, also on the chain level.

\begin{df}
\label{quasiopdef}
A quasi-operad is an operad in which the axiom of associativity
need not hold and a quasi-PROP is a PROP is which the axiom of
associativity need not hold.

A quasi-operad/PROP is  called a homotopy operad/PROP, if it is in
the category of topological spaces and the associativity equations
(respectively the compatibility equations) hold up to homotopy.
\end{df}

\subsubsection{The Sullivan quasi-PROP}

Notice that for an element $\a \in \iooarc$, we can independently
scale the weights of the arcs running to the inputs, since there
are no arcs which go from one input to another input. Therefore
the set of all arcs are in a 1-1 correspondence to the union of
the sets of arcs incident to the $\In$ boundaries. In order to glue
the inputs to outputs, using this bijection we first scale the
weights on the input factors locally
---that is   separately for each boundary homogeneously  we scale all the
weights of edges incident to that boundary--- to make them match
with the output weights and then glue. More precisely, let $\a\in
\iooarc(n,m)$  and $\b \in \iooarc(m,k)$.
 Let $\Out(\a)$ be the indexing
set of the ``out'' boundaries of $\a$ and $\In(\b)$ be the indexing set
of the ``in'' boundaries of $\b$. For $j\in \Out(\a)$ and $i\in
\In(\b)$
 let $w_j$  be the sum of the weights incident
to the boundary $j$ of $\a$ and similarly let $w_i$ be the sum of
the weights of the arcs incident to the boundary $i$ of $\b$. Let
$\phi: \Out(\a) \rightarrow \In(\b)$ be the bijection for the
gluing. Then we define $\a \bullet_{\phi}\b$ to be the arc family
obtained by first scaling the arcs incident
 to the boundary
$i$ of $\b$ by the factor $\frac{w_j}{w_i}$ (where
$j=\phi^{-1}(i)$), and then gluing the arcs of $\a$ and $\b$ along
the boundaries $i$ and $j$ which now have the same total incident
weight.

To define the vertical compositions one can use disjoint union of
surfaces by passing to possibly disconnected surfaces.

\begin{thm}
The operations $\bullet_{\phi}$  imbue the sets generated
by $\Sn
\times \Sm$ modules $\iooarc(n,m)$
with the structure of a homotopy PROP.
\end{thm}

\begin{proof}
Tedious, but straightforward. The main observation is that
if one lets $C_i$ be the images of the boundaries
after the gluing of the surfaces, each of the glued arcs can be
moved to transversally cut exactly one of these curves. Now
there are commuting flows, that act locally at $C_i$ which
essentially move the weights of the arcs before gluing. The example
of one boundary gluing is given in \cite{KP}. The flow is a
generalization of that of \cite{cact} and  given by
moving from one partial homotopy diagonal to another.
These flows allow one to flow from one association
to the other.
 \end{proof}

\begin{cor}
The homology of  $\iooarc(n,m)$ is a PROP.
\end{cor}

\begin{nota}
We phrased the theorem so that we can
avoid yet other notation. What is meant is
that the operations $\bullet_{\phi}$ and the $\Sn\times\Sm$
actions naturally induce operations and permutation actions
 on the PROP whose $(N,M)$--component is given by
$$\coprod_{(n_1,\dots, n_k), \sum n_i=N; (m_1,\dots m_k),\sum n_j=M}
\prod_{n_i, m_j} \iooarc(n_i,m_i)$$ We will call this PROP simply
$\iooarc$ and also use this shorthand notation for any sub-PROP
generated by a subset of $\iooarc$.
\end{nota}

\begin{rmk}
This scaling when restricted to $\Tree$ or $\Lintree$ is the
quasi--operad structure used for the normalized versions of cacti
and spineless cacti \cite{cact}.
\end{rmk}
\begin{rmk}
\label{connectedremark}
If one wishes to avoid disconnected surfaces, one can use the
following trick. Let $I(m,m)\in \ioarc(m,m)$ be the surface of
genus 0 with $2m$ boundaries which are partitioned into two sets
$\In$ and $\Out$ which are individually numbered from $1$ to $m$
with exactly one arc running from the in
 boundary $i$ to the out boundary $i$. In this case,
  we define the vertical PROP
composition $\a\otimes \b$ for $\a\in \iooarc(m,n)$ and $\b\in
\iooarc(k,l)$
 to be given by the gluing of the out boundaries of $\a$
onto the first $n$ $\In$ boundaries of $I(n+l,n+l)$  and gluing the
$\Out$ boundaries of $\b$ to the last $l$ of the $\In$ boundaries of
$I(n+l,n+l)$.
\end{rmk}

\subsection{Cell models}
Without reiterating all the details of the construction, it is
clear that there are  induced (quasi)--operads, PROPs and
di-operads on the graded open cell level. These structures can be
seen to be strict on the this level. If we pass to the graded
cells also the rational operad structure becomes a strict operad structure.
We keep the convention in the PROP gluing that
is analogous to equation \ref{setop}, namely that the operation
on cells is zero if for any of the pairs of boundaries which are glued
in the PROP operation both boundaries are twisted.

\begin{thm}
\label{cellthm}
On the (graded) open cell level the homotopy (or rational)
structures give rise to
the non-homotopy structures, moreover the quasi--filling subspaces
are isomorphic to types of graphs and hence there is the induced
structure on these graphs. In particular,\\

\begin{tabular}{l|l}
&is a\\
\hline
$\OC(\Diarc)$&di--operad\\
$\OC(\iooarc)$&sub--di--operad and a PROP\\
$\OC(\ioarc)$&sub--di--operad and a sub--PROP\\[2mm]
\end{tabular}\\
and furthermore going to the graded cells, we obtain operads
which are isomorphic to the free Abelian groups of certain
types of ribbon graphs\\[2mm]

\begin{tabular}{l|l|l}
& is a& and is isomorphic to\\
\hline
$\Gr\OC(\Anarc)$&cyclic operad\\
$\Gr\OC(\Anarcn)$&cyclic operad&$\Anrib$\\
$\Gr\OC(\Arcn)$&dg cyclic operad&$\Rib$\\
$\Gr\OC(\ioarc_{\#})$&dg PROP&$\ppdigraph$\\[2mm]
\end{tabular}

Here the table is to be understood in the sense that the entry in
the first column carries the structure of the second column and is
additively (over $\Z$) isomorphic to the third column thereby
inducing the respective structure.
\end{thm}

\begin{proof}
The claims for the di--operad structures follow for the results about $\Arc$.
For the PROP
structure, the arguments are analogous to those of \S \ref{cellop}.
 We again claim that the PROP operations are 1--1
on cells which are not twisted at the boundary. The case of closed
loops cannot appear. After scaling, the gluing operation is locally
given by the shuffle combinatorics of \S\ref{gluingpar}, so that
indeed the image of the PROP action will be full cells and 1-1 in
the case that not both the boundaries are twisted. This local
argument accounts for all arcs of the glued surface passing through
the images of the curves defined by the glued boundaries whose
weights can all be scaled independently, while all other arcs are
unchanged. So indeed we get an induced map on cells. Locally, that
is regarding each pair of glued boundaries separately, the map is
1--1 precisely if the two boundaries are not both twisted and there
are no closed loops. This yields the above assertion. The fact that
the image cells lie in the relevant PROP is clear by the definition
of the restriction. Combining two ``in'' to ``out'' arcs yields an
arc of the same type, and also combining an ``out'' to ``out'' arc
with an ``in'' to ``out'' arc again yields an arc running from
``out'' to ``out''. For the rational structures, we notice that
again all ``problems'' arise in at most codimension one.
Consequently  passing to the associated graded induces the
unrestricted structures.
\end{proof}

\subsection{A CW model for Sullivan Chord diagrams}
Recall that $\A$ was the CW complex that contains $\Arc$ as
a subset. We define  $\A^{s\; i/o}_{g,r}:=\A^s_{g,r}\times
(\Z/2\Z)^r$ analogously to $\Diarc$. Notice that these spaces are
CW-complexes. We call the collection of these spaces $\DiA$ and as
usual write $\DiA(n):=\amalg_{g,s} A^{s\; i/o}_{g,n+1}$.

\begin{df}
We let $\Diiooarc$ be the collection of subspaces of the spaces of
$\DiA$ in which there are only arcs running from the $\In$ to
the $\Out$ and possibly
from the $\Out$ to the $\Out$
boundaries  {\em and} there is no empty $\In$ boundary.

We let $\Diioarc\subset\Diiooarc$ be the subset where arcs only run
from the $\In$ boundaries to the $\In$ boundaries.

We also let $\LDiiooarc$ be the subset of elements whose
underlying arc graph is not twisted at any of the $\In$ boundaries
and set $\LDiioarc=\LDiiooarc\cap\Diioarc$.
\end{df}

\begin{rmk}
It is clear that $\Diioarc$ is a retract of $\Diiooarc$ by simply
scaling all the weights on the arcs connecting $\Out$ to $\Out$ to zero.
So any cell model for $\Diioarc$ also yields a cell model for $\Diiooarc$.
\end{rmk}

\subsubsection{Extended gluing on $\Diioarc$}
\label{dgextend} We can extend the operadic gluing to $\DiA$ as in
\S\ref{extdsec}, in particular for $\Diiooarc$ this means that if
we glue to an empty ``out'' boundary, the arcs incident to the
respective ``in''  boundary will be deleted. Notice that this
leaves us in $\Diiooarc$ since none of the arcs hitting the ``in''
boundaries in the glued surface have been effected. As mentioned
earlier this operation is only associative up to homotopy, so
$\Diiooarc$ is a quasi--di--operad. With respect to the PROP
structure we define the maps $\bullet_{\phi}$ by setting the
weights of the arcs incident to an ``in'' boundary that is glued
to an empty ``out'' boundary to zero. These gluings then goes over
to the cell level as in \S\ref{cellops} and Theorem \ref{cellthm}.

In particular, this yields the following extension  of the gluing to
 $\OC(\Diiooarc)$. The gluing of
a cell indexed by an arc graph with an empty ``out''
glued at that ``out'' to a non-empty
``in'' is defined to be the cell indexed by the modified image of the two
arc graphs, where the modification is that
in the second arc graph all edges incident to
the boundary have been deleted as detailed in \ref{extasec}.

Just like for cacti and spineless cacti \cite{cact}, there is a
smaller space which is a retract of $\Diioarc$ that is actually a
CW complex.

\begin{df}
We define $\Diioarci\subset \Diioarci$ to be the subspace of graphs
whose sum of weights of arcs incident to every $\In$ boundary
vertex is one.
\end{df}

It is clear that $\Diioarc$ retracts to $\Diioarci$ by
homogenously scaling the weights of the sets  of arcs incident to
each $\In$ boundary separately for each of these sets, so that
their total weight becomes one.

\begin{prop}
$\Diioarci$ is a CW complex, whose cells are indexed by the arc
graphs of the given type.
\end{prop}
\begin{proof}
Completely analogous to the constructions of \cite{cact,del}.
Given an arc graph, we define the cell by $\times_{v\in \In}
\Delta^{|v|}$ where $\D$ is the standard simplex. We define the
attaching maps by gluing the boundary corresponding to a
face of a simplex to the cell of lower dimension indexed by
the arc graph obtained form the original arc graph by
deleting exactly the edge that indexes the face. It is then straightforward
to show that this CW complex realizes $\Diioarci$.
\end{proof}

It is clear that the graphs of $\LDiioarc$ form a sub--CW complex
of $\Diioarci$ which we call $\LDiioarci$.

We define the di--operadic  compositions  on $\Diiooarci$ by scaling
the input individually to the weight of the output. The homotopy PROP structure
is just  the homotopy sub--PROP structure. Notice that in the gluings
$\bullet_{\phi}$ one only scales at the $\In$ boundaries which are
to be glued so that the
weights on the $\In$ boundaries which remain after gluing are unchanged.

\begin{thm}
These compositions define a homotopy--PROP structure on $\Diiooarci$
which makes $\Diioarc$, $\LDiioarc$ and $\LDiiooarc$ into a
sub--quasi--PROPs. Their homotopy--PROP structures descend to PROP
structures on $\OC(\Diiooarc)$ and $\Gr\OC(\Diiooarc)$. Also, the
PROP structure of $\Diioarci$
 descends to a $dg$--PROP structure on $CC_*(\Diioarci)$.
Moreover, the two PROP structures $CC_*(\Diioarci)$ and $\Gr\OC(\Diioarc)$
and  their
differentials when viewed as defined on the same free Abelian group generated
by the respective arc graphs  agree.

The same statement {\em mutatis mutandis} holds true for the
respective di--operad structures. Likewise the analogous statements
hold true for $\LDiioarc$ and $\LDiioarci$ as well.
\end{thm}

\begin{proof}
Tedious but straightforward generalization  to the
case of several outputs
 of the analogous statement about $\Cacti^1$ and $\Cacti$ given in \cite{cact}.
Since we are dealing with a local scaling which acts independently
on the arcs due to the restrictions we imposed, one can construct a
homotopy which scales the sum of the weights of the arcs through any
given fixed closed curve which is the image of a boundary under a
gluing to say $1$. Using these homotopies, one can flow from one
association to the other, this shows the quasi--PROP structure. On
the cell level, one has to make sure that all the combinatorially
possible graphs arise and that each weight only arises once. This is
a straightforward verification using the techniques presented above.
The sub--PROPs are actually stable under gluing, since one will
never induce a twist on the $\In$ boundaries when gluing. Likewise
one cannot obtain any ``out'' to ``out'' arcs if they were not
previously there. The last statement about the CW--complex follows
from the fact that the cells of lower dimension are killed in the
cellular chain complex.
\end{proof}

\begin{cor}
The $\Sn\times \Sm$--modules  $CC_*(\Diioarci)(n,m)$ form a
$dg$--PROP and give a chain model operad for $\Diioarc$ that is
for the  extended metric Sullivan Chord diagrams.

The same statement holds true for the respective di--operad structures.

\end{cor}\qed

\begin{cor}
$H_*(\Diioarci)\simeq H_*(\Diioarc)$ and
the induced  PROP structures agree. Hence $\Diioarci$ is a cellular--PROP
 model for the extended Sullivan Chord PROP $\Gr\OC(\Diioarc)$.

The same statement holds true for the respective di--operad
structures.
\end{cor}\qed


\begin{thebibliography}{99}


\bibitem
[BF] {BF} C.\ Berger and B.\ Fresse {\it Une d\'{e}composition
prismatique de l'op\'{e}rade de Barratt-Eccles.}  C. R. Math. Acad.
Sci. Paris 335 (2002), no. 4, 365-370.



\bibitem
[C]  {C}
R.~ Cohen. {\it Multiplicative properties of Atiyah duality.}  Homology Homotopy Appl.  6  (2004),  no. 1, 269--281
\bibitem
[CG]  {CG}
R.~L.~Cohen and V.~Godin.{\it A polarized view of string topology.}  Topology, geometry and quantum field theory,  127--154, London Math. Soc. Lecture Note Ser., 308, Cambridge Univ. Press, Cambridge, 2004.
\bibitem
[Ch]  {Chat}
D.~Chataur. {\it A bordism approach to string topology.}  Int. Math. Res. Not.  2005,  no. 46, 2829--2875.

\bibitem
[CJ]  {CJ}
R.\ L.\ Cohen and J.D.S.\ Jones
{\it  A homotopy theoretic realization of string topology}
Math.\ Ann.\ 324 (2002), no. 4, 773-798.
 \bibitem
[CS]  {CS} M.\ Chas and D.\ Sullivan. {\it  String Topology.}
Preprint math.GT/9911159. To appear in Ann.\ of Math.\

\bibitem
[Co]  {cost}
 K.~J.~Costello.
{\it Topological conformal field theories and Calabi-Yau categories}.
Preprint math.QA/0412149.
 {\it The Gromov-Witten potential associated to a TCFT.}
Preprint math.QA/0509264 and {\it A dual point of view on the ribbon graph decomposition of moduli space.}
Preprint math.GT/0601130.
\bibitem
[CV]{Vogt}
J.~Conant and K.~Vogtmann. {\it On a theorem of Kontsevich.}
 Algebr.\ Geom.\ Topol.\  3  (2003), 1167--1224


\bibitem
[H]{Harstab}
J~.L Harer.{\it Stability of the homology of the mapping class
groups of orientable surfaces},  { Ann. of Math.}  {  121}
(1985), 215-249.

\bibitem
[HM] {HubMas}
J.~H. Hubbard and H.~Masur, {\it Quadratic differentials and
foliations},  {Acta Math.}  {142} (1979), 221-274.


\bibitem
[J]{jones} J.D.S.~Jones. {\it Cyclic homology and equivariant homology.}
 Inventionnes Math. 87 (1987), 403--423

\bibitem
[K1]{cact} R.\ M.\ Kaufmann. {\it On several varieties of cacti
and their relations.} Algebraic \& Geometric Topology 5 (2005),
2--300.

\bibitem
[K2]{del} R.\ M.\ Kaufmann. {\it  On Spineless Cacti, Deligne's
Conjecture and Connes--Kreimer's Hopf Algebra. } Preprint
math.QA/0308005.


\bibitem
[K3]{cyclic}
R.~M.~Kaufmann.
{\it A proof of a cyclic version of Deligne's conjecture via Cacti.} Preprint,  math.QA/0403340.

\bibitem
[K4]{hoch2} R.\ M.\ Kaufmann. {\it Moduli space actions on the Hochschild Co-Chains of a Frobenius
algebra II: Correlators}. Preprint.


\bibitem
[K5]{Ribbon} R.\ M.\ Kaufmann. {\it The Arc Spectrum.} In
preparation.

\bibitem
[KLi1]{KLi1}
A.~Kapustin and Y.~Li, {\em D-branes in Landau-Ginzburg models and
algebraic geometry}. JHEP {\bf 0312}, 005 (2003)

\bibitem
[KLi2]{KLi2}
A.~Kapustin and Y.~Li, {\em Topological correlators in
Landau-Ginzburg models with boundaries}. Preprint hep-th/0305136.




\bibitem
[KR]{KR}  A.~Kapustin and L.~Rozansky.
{\em  On the relation between open and closed topological
strings.} Commun.Math.Phys. 252 (2004) 393-414

\bibitem
[Ko1]{Koint}
M.~Kontsevich.
{\it Intersection theory on the moduli space of curves and the matrix Airy function.}  Comm. Math. Phys.  147  (1992),  no. 1, 1--23

\bibitem
[Ko2]{formal} M.~Kontsevich. {\it Formal (non)commutative
symplectic geometry.} The Gel'fand Mathematical Seminars,
1990--1992,  173--187, Birkh\"auser Boston, Boston, MA, 1993.

\bibitem
[Ko3]{Maxim}
M.\ Kontsevich.
{\it Operads and Motives in Deformation Quantization.}
Lett.Math.Phys.\ 48 (1999) 35-72.


\bibitem
[KP]{KP} R.\ M.\ Kaufmann and R.\ B.\ Penner. {\it
Closed/open string diagrammatics}. math.GT/0603485.
To appear in Nucl.~Phys.~B.

\bibitem
[KLP]{KLP} R.\ M.\ Kaufmann, M.\ Livernet and R.\ B.\ Penner. {\it
Arc Operads and Arc Algebras.} Geometry and Topology 7 (2003),
511-568.


\bibitem
[KS1]  {KS} M.\ Kontsevich and Y.\ Soibelman. {\it Deformations of
algebras over operads and Deligne's conjecture.} Conf\'erence
Mosh\'e Flato
  1999, Vol. I (Dijon), 255--307,
Math. Phys. Stud., 21, Kluwer Acad. Publ., Dordrecht, 2000.

\bibitem
[KS2]  {KS2} M.\ Kontsevich and Y.\ Soibelman.
 {\it Notes on $\Ainf$--categories
and non--commutative geometry}. Manuscript.




\bibitem
[Me]{Merk}
S.~A.~Merkulov.{\it  De Rham model for string topology.}  Int.\ Math.\
Res.\ Not.\  2004,  no.\ 55, 2955--2981.


\bibitem
[MS1]  {MS}
J.\ E.\ McClure and J.\ H.\ Smith, Jeffrey H.
{\it  A solution of Deligne's Hochschild cohomology conjecture.}
Recent progress in homotopy
  theory (Baltimore, MD, 2000), 153-193,
Contemp. Math., 293, Amer. Math. Soc., Providence, RI, 2002.

\bibitem
[MS2]{MS2}
J.\ E.\ McClure and J.\ H.\ Smith, Jeffrey H.
{\it Multivariable cochain operations and little $n$-cubes.}
J. Amer. Math. Soc. 16 (2003), no. 3, 681--704


\bibitem
[MS3]  {MScosimp} James E.~McClure and Jeffrey H.~Smith.
{\it Cosimplicial objects and little $n$-cubes. I.}  Amer. J. Math.  126  (2004),  no. 5, 1109--1153.


\bibitem
[P]{P3} R.\ C.\ Penner, ``The decorated Teichm\"uller space of
punctured surfaces", {\it Communications in Mathematical  Physics}
{\bf 113}   (1987),  299-339.


\bibitem
[P2] {P2} R.\ C.\ Penner. {\it Decorated Teichm\"uller theory of
bordered surfaces.}
 Comm. Anal. Geom.  12  (2004),  no. 4, 793--820.

\bibitem
[S1]{C1}
D.~Sullivan. {\em Sigma models and string topology.} Graphs and
patterns in mathematics and theoretical physics,  1--11, Proc.
Sympos. Pure Math., 73, Amer. Math. Soc., Providence, RI, 2005.
(Reviewer: David Chataur)
\bibitem
[S2]{C2}
 D.~Sullivan. {\em Open and closed string field theory interpreted
in classical algebraic topology.}  Topology, geometry and quantum
field theory,  344--357, London Math. Soc. Lecture Note Ser., 308,
Cambridge Univ. Press, Cambridge, 2004



\bibitem
[St]
{Strebel}
K. Strebel, {\it Quadratic Differentials}, { Ergebnisse der Math.}
{ 3:5}, Springer-Verlag, Heidelberg (1984).


\bibitem
[T] {T} Tamarkin, D., {\it Another proof of M. Kontsevich formality
theorem}. Peprint math/9803025.\\
{\it Formality of Chain Operad of Small Squares}. Lett. Math. Phys. 66 (2003), no. 1-2, 65--72.

\bibitem
[TZ] {TZ}
T.~Tradler and M.~Zeinalian.
{\it On the cyclic Deligne conjecture}
J. Pure Appl. Algebra 204 (2006), no. 2, 280--299.
\bibitem
[V]  {Vor2}
A.\ A.\  Voronov.
{\it  Homotopy Gerstenhaber algebras. }
Conf\'erence Mosh\'e Flato 1999, Vol. II (Dijon), 307-331, Math. Phys.
  Stud., 22, Kluwer Acad. Publ., Dordrecht, 2000
\bibitem
[V2]  {Vor} A.\ A.\  Voronov. {\it  Notes on universal algebra.}
Graphs and Patterns in Mathematics and Theoretical Physics (M.
Lyubich and L. Takhtajan, eds.), Proc. Sympos. Pure Math., vol.
73. AMS, Providence, RI, 2005, pp. 81-103.

\end{thebibliography}
\end{document}